\documentclass[a4paper,12pt]{amsart}
\pdfoutput=1
\usepackage{setspace}
\usepackage[left=35mm,head=30mm,bottom=30mm,foot=20mm,right=35mm]{geometry}
\usepackage{amsmath,amsthm,amsfonts,amssymb,mathtools,mathrsfs,url,bm,enumitem,stackengine}
\usepackage{newtxtext}
\usepackage[varvw]{newtxmath}
\allowdisplaybreaks
\usepackage{hyperref}
\hypersetup{
  colorlinks=true,
  linkcolor=blue,
  citecolor=red,
}
\mathtoolsset{showonlyrefs=true}
\theoremstyle{definition}
  
  \newtheorem{rem}{Remark}[section]
\theoremstyle{plain}
  \newtheorem{thm}{Theorem}[section]
  \newtheorem{cor}{Corollary}[section]
  \newtheorem{prop}{Proposition}[section]
  \newtheorem{lem}{Lemma}[section]

\renewcommand{\epsilon}{\varepsilon}
\newcommand\restr[2]{{
  \left.\kern-\nulldelimiterspace 
  #1 
  \vphantom{\big|} 
  \right|_{#2} 
}}

\begin{document}

\title[Time-asymptotic expansion]{Time-asymptotic expansion with pointwise remainder estimates for 1D viscous compressible flow}
\author{Kai Koike}
\address{Department of Mathematics, Tokyo Institute of Technology, Tokyo 152-8551, Japan}
\email{koike.k@math.titech.ac.jp}
\date{\today}

\begin{abstract}
  We construct a time-asymptotic expansion with pointwise remainder estimates for solutions to 1D compressible Navier--Stokes equations. The leading-order term is the well-known diffusion wave and the higher-order terms are newly introduced family of waves which we call \textit{higher-order diffusion waves}. In particular, these provide accurate description of the power-law asymptotics of the solution around the origin $x=0$ where the diffusion wave decays exponentially. The expansion is valid locally and also globally in the $L^p(\mathbb{R})$-norm for all $1\leq p\leq \infty$. The proof is based on pointwise estimates of Green's function.
\end{abstract}

\maketitle


\section{Introduction}
The equations
\begin{equation}
  \label{eq:p_system}
  \begin{dcases}
    v_t-u_x=0,                      & x\in \mathbb{R},\, t>0, \\
    u_t+p(v)_x=\nu(u_x/v)_x,        & x\in \mathbb{R},\, t>0, \\
    v(x,0)=v_0(x),\, u(x,0)=u_0(x), & x\in \mathbb{R}
  \end{dcases}
\end{equation}
describe the motion of a 1D viscous compressible flow. Here $v(x,t)$ is the specific volume (the reciprocal of the density $\rho$) and $u(x,t)$ is the flow velocity; $t$ is the time and $x$ is the Lagrangian mass coordinate related to the Eulerian coordinate $X$ by $x=\int_{0}^{X}\rho(X',t)\, dX'$. The system above models barotropic flow, that is, the pressure $p(v)$ does not depend on the temperature. We assume that $p'(v)<0$ and $p''(v)\neq 0$ for $v>0$ and that the viscous coefficient $\nu$ is a positive constant. The system is often called the $p$-system in the literature and is a typical example of quasilinear hyperbolic-parabolic viscous conservation laws.

The purpose of this paper is to construct a time-asymptotic expansion of the solution to~\eqref{eq:p_system} together with pointwise estimates for the remainder. We shall consider solutions close to the steady state $(v_S,u_S)\equiv (1,0)$. To study the long-time asymptotics of such solutions, it is convenient to consider
\begin{equation}
  \label{eq:ui}
  u_1=\frac{p''(1)}{4c}[-(v-1)+u/c], \quad u_2=\frac{p''(1)}{4c}[(v-1)+u/c]
\end{equation}
instead of $(v,u)$. Here $c=\sqrt{-p'(1)}$ is the speed of sound for the state $(v_S,u_S)$.

It is well-known that $u_i$ has the \textit{diffusion wave} $\theta_i$ as its asymptotic profile. Here $\theta_i$ ($i=1,2$) is the self-similar solution to the convective viscous Burgers equation
\begin{equation}
  \label{eq:thetai}
  \begin{dcases}
    \partial_t \theta_i+\lambda_i \partial_x \theta_i+\partial_x(\theta_{i}^{2}/2)=\frac{\nu}{2}\partial_{x}^{2}\theta_i, & x\in \mathbb{R},t>0, \\
    \lim_{t\searrow -1}\theta_i(x,t)=M_i \delta(x), \quad x\in \mathbb{R},
  \end{dcases}
\end{equation}
where $\lambda_i=(-1)^{i-1}c$, $M_i=\int_{-\infty}^{\infty}u_i(x,0)\, dx$, and $\delta(x)$ is the Dirac delta function. An explicit formula for $\theta_i$ is available through the use of Cole--Hopf transformation:
\begin{equation}
  \label{eq:theta_explicit}
  \theta_i(x,t)=\frac{\sqrt{\nu}\left( e^{\frac{M_i}{\nu}}-1 \right)}{\sqrt{2(t+1)}}e^{-\frac{(x-\lambda_i(t+1))^2}{2\nu(t+1)}}\left[ \sqrt{\pi}+\left( e^{\frac{M_i}{\nu}}-1 \right) \int_{\frac{x-\lambda_i(t+1)}{\sqrt{2\nu(t+1)}}}^{\infty}e^{-y^2}\, dy \right]^{-1}.
\end{equation}
The diffusion wave $\theta_i$ describes the leading-order asymptotics in the $L^p(\mathbb{R})$-norm. In fact, we have the following optimal decay estimates~\cite{Zeng94}:
\begin{equation}
  \| \theta_i(\cdot,t) \|_{L^p}\lesssim t^{-(1-1/p)/2} \quad \text{and} \quad \| (u_i-\theta_i)(\cdot,t) \|_{L^p}\lesssim t^{-(3/2-1/p)/2} \quad (1\leq p\leq \infty).
\end{equation}
The key to prove these $L^p$-decay estimates---especially for $p=1$---is the pointwise estimates for Green's function of the linearization of~\eqref{eq:p_system} around $(v_S,u_S)$. They in fact allow us to obtain pointwise estimates for the solution itself~\cite{LZ97}:
\begin{equation}
  \label{eq:ui-thetai_PWE}
  |(u_i-\theta_i)(x,t)|\lesssim  [(x-\lambda_i(t+1))^2+(t+1)]^{-3/4}+[|x+\lambda_i(t+1)|^3+(t+1)^2]^{-1/2}.
\end{equation}
The $L^p$-decay estimates are obtained by integrating this.

Pointwise estimates~\eqref{eq:ui-thetai_PWE} allow us to deduce not just global $L^p$-estimates but also local ones. In particular, we have $|(u_i-\theta_i)(x,t)|\lesssim t^{-3/4}$ for $x=\lambda_i t+O(1)$. Since $\theta_i(x,t)\lesssim t^{-1/2}$ for $x=\lambda_i t+O(1)$, the diffusion wave $\theta_i$ also describes the leading-order asymptotics locally around the characteristic line $x=\lambda_i t$. However, the situation is different around the origin $x=0$. As can be seen from~\eqref{eq:theta_explicit}, the diffusion wave $\theta_i$ decays exponentially fast around the origin $x=0$ but~\eqref{eq:ui-thetai_PWE} implies $|(u_i-\theta_i)(x,t)|\lesssim t^{-3/2}$ for $x=O(1)$. Thus the diffusion wave $\theta_i$ provides almost no information about the long-time asymptotics around $x=0$; we need new waves to capture the asymptotic behavior there.

In~\cite{vBPW08}, van Baalen, Popovi\'{c}, and Wayne constructed a time-asymptotic expansion of $u_i$ in an $L^2$-framework. The leading-order term of the expansion is the diffusion wave $\theta_i$ but the first higher-order term beyond $\theta_i$ turns out to be a wave decaying algebraically as $t^{-3/2}$ around the origin. It is then natural to expect that this new wave captures the leading-order asymptotics of the flow around $x=0$. However, the decay estimate for the remainder of the expansion is given in the $H^1(\mathbb{R})$-norm. This implies only a far from optimal decay estimate around $x=0$. For this reason, we cannot conclude that the higher-order term describes the leading-order asymptotics of the flow around $x=0$.

To overcome this issue, we construct a time-asymptotic expansion of $u_i$ with pointwise estimates for the remainder. The leading-order term is the diffusion wave $\theta_i$ and the higher-order terms are \textit{higher-order diffusion waves} $\xi_{i;n}$ ($n\geq 1$) defined in the next section. It turns out that $|\xi_{i;n}(x,t)|\lesssim t^{-(2-1/2^n)}$ for $x=O(1)$ as $t\to \infty$. Setting $n=1$, we see that $|\xi_{i;1}(x,t)|\lesssim t^{-3/2}$ for $x=O(1)$. The pointwise estimates for the remainder imply $|(u_i-\xi_{i;1})(x,t)|\lesssim t^{-7/4}$ for $x=O(1)$, thus it is rigorously proved that $\xi_{i;1}$ describes the leading-order asymptotics of $u_i$ for $x=O(1)$. In addition, thanks to the pointwise estimates, our asymptotic expansion is valid not only in the $L^2(\mathbb{R})$-norm but also in the $L^1(\mathbb{R})$-norm.

The proof is based on pointwise estimates of Green's function, and the basic strategy follows that of~\cite{LZ97}. The most non-trivial part of the proof is perhaps the definition of the higher-order diffusion waves $\xi_{i;n}$ ($n\geq 1$); see~\eqref{eq:xii_n}. Although the differential equation defining $\xi_{i;n}$ does not seem to have a simple solution formula such as~\eqref{eq:theta_explicit},\footnote{Nevertheless, we provide accurate asymptotic analysis in Propositions~\ref{prop:xii_n} and~\ref{prop:xii_n:fn}.}~we use its structure (by the help of Lemma~\ref{lem:heat_equation_technique}) to analyze cancellation effects which are crucial in nonlinear estimates; see the proof of Lemma~\ref{lem:PWE_nonlinear}.

Before closing the introduction, we briefly comment on related works. Diffusion wave approximations and pointwise estimates of solutions has been extensively studied for hyperbolic-parabolic systems~\cite{LZ97}, hyperbolic-elliptic systems~\cite{IK02}, hyperbolic balance laws~\cite{WY03,ZC16}, the Boltzmann equation~\cite{LY04}, and so on. In these works, nonlinear diffusion waves similar to $\theta_i$ were constructed and pointwise estimates of solutions were obtained. But to the author's knowledge, time-asymptotic expansions with pointwise estimates have not been obtained previously. We mention that the author already analyzed the second-order term $\xi_{i;1}$ in connection with a fluid--structure interaction problem in~\cite{Koike22b}; the complete asymptotic expansion, however, was not given. We also comment that for multidimensional incompressible Navier--Stokes equations, time-asymptotic expansions were studied for example in~\cite{Carpio96,FM01}. Because the nonlinearity is weaker compared to the 1D case, nonlinear waves similar to $\xi_{i;n}$ do not appear in these works.

In the next section, we state our main results. These are proved in Section~\ref{sec:proof}.

\section{Main results}
To state our main results (Theorem~\ref{thm:main1}) we start by defining and discussing the properties of the higher-order diffusion waves $\xi_{i;n}$ ($n\geq 1$) mentioned in the introduction.

\subsection{Higher-order diffusion waves}
\label{sec:xii_n}
Let $(v,u)$ be the solution to~\eqref{eq:p_system}. Then define $u_i$ by~\eqref{eq:ui} and set
\begin{equation}
  \label{eq:Mi}
  M_i=\int_{-\infty}^{\infty}u_i(x,0)\, dx.
\end{equation}
Let $\xi_{i;0}=\theta_i/2$ with $\theta_i$ defined by~\eqref{eq:thetai}. We then define the higher-order diffusion waves $\xi_{i;n}$ ($n\geq 1$) inductively by the equations
\begin{equation}
  \label{eq:xii_n}
  \begin{dcases}
    \partial_t \xi_{i;n}+\lambda_i \partial_x \xi_{i;n}+\partial_x(\theta_i \xi_{i;n})+\partial_x(\theta_{i'}\xi_{i';n-1})=\frac{\nu}{2}\partial_{x}^{2}\xi_{i;n}, & x\in \mathbb{R},t>0, \\
    \xi_{i;n}(x,0)=0, & x\in \mathbb{R}.
  \end{dcases}
\end{equation}
Here $\lambda_i=(-1)^{i-1}c$ and $i'=3-i$, that is, $1'=2$ and $2'=1$. We remind the reader that $c=\sqrt{-p'(1)}>0$.

Although we do not have a simple explicit formula for $\xi_{i;n}$ like~\eqref{eq:theta_explicit}, we can still understand its asymptotic behavior quite well. To explain this, we introduce
\begin{equation}
  \label{eq:alpha_beta}
  \alpha_n=2-\frac{1}{2^{n+1}}, \quad \beta_n=\frac{3}{2}-\frac{1}{2^{n+1}} \quad (n\geq -1)
\end{equation}
and
\begin{align}
  \label{eq:aux}
  \begin{aligned}
    \psi_n(x,t;\lambda) & =[(x-\lambda(t+1))^2+(t+1)]^{-\alpha_n/2}, \\
    \tilde{\psi}_{n}(x,t;\lambda) & =[|x-\lambda(t+1)|^{\alpha_n}+(t+1)^{\beta_n}]^{-1}, \\
    \Psi_{i;n}(x,t) & =\psi_n(x,t;\lambda_i)+\tilde{\psi}_n(x,t;\lambda_{i'}).
  \end{aligned}
\end{align}
Then we have the following decay estimates for $\xi_{i;n}$ (we postpone the proof until we later prove a finer version in Lemma~\ref{lem:xii_n}).

\begin{prop}
  \label{prop:xii_n}
  Let $n\geq 1$ and $\epsilon=\max(M_1,M_2)$. For $k\geq 0$, if $\epsilon$ is sufficiently small, we have
  \begin{equation}
    |\partial_{x}^{k}\xi_{i;n}(x,t)|\leq C_{n,k}\epsilon^{n+1}(t+1)^{-k/2}\Psi_{i;n-1}(x,t)
  \end{equation}
  for some positive constant $C_{n,k}$. In particular, when $|x|\leq K$ for some fixed $K>0$, we have
  \begin{equation}
    |\xi_{i;n}(x,t)|\leq C_{n,K}(t+1)^{-\alpha_{n-1}}
  \end{equation}
  for some $C_{n,K}>0$. Moreover, for any $1\leq p\leq \infty$, there exists $C_{n,p}>0$ such that
  \begin{equation}
    \| \xi_{i;n}(\cdot,t) \|_{L^p}\leq C_{n,p}(t+1)^{-(\alpha_{n-1}-1/p)/2}.
  \end{equation}
\end{prop}

We can also prove more detailed estimates if we focus on $x$ with $(-1)^{i-1}x\geq 0$. Let
\begin{equation}
  \label{eq:g}
  g(z)=\partial_x e^{-\frac{x^2}{2\nu}}=-(x/\nu)e^{-\frac{x^2}{2\nu}},
\end{equation}
\begin{equation}
  \label{eq:f0}
  f_{i;0}(z)=\frac{\sqrt{\nu}}{\sqrt{2}}\left( e^{\frac{M_i}{\nu}}-1 \right) e^{-\frac{z^2}{2\nu}}\left[ \sqrt{\pi}+\left( e^{\frac{M_i}{\nu}}-1 \right) \int_{z}^{\infty}e^{-\xi^2}\, d\xi \right]^{-1},
\end{equation}
and
\begin{equation}
  \label{eq:fn}
  f_{i;n}(z)=\int_{(-1)^{i-1}z}^{\infty}[\xi-(-1)^{i-1}z]^{-(1-1/2^n)}\xi e^{-\frac{\xi^2}{2\nu}}\, d\xi.
\end{equation}
We then have the following asymptotic formula. This is obtained from Lemma~\ref{lem:xi_n:etai_n} proved in the next section.

\begin{prop}
  \label{prop:xii_n:fn}
  Let $n\geq 1$ and $\epsilon=\max(M_1,M_2)$. For any $K>0$, if $\epsilon$ is sufficiently small, there exist $A_{i;n}$, $B_{i;n}$, and $C_n>0$ such that
  \begin{align}
    & \left| \left\{ \xi_{i;n}(x,t)-\frac{A_{i;n}}{(t+1)^{\alpha_{n-1}/2}}f_{i;n}\left( \frac{x-\lambda_i(t+1)}{\sqrt{t+1}} \right)-\frac{B_{i;n}}{(t+1)^{\alpha_{n-1}/2}}g\left( \frac{x-\lambda_i(t+1)}{\sqrt{t+1}} \right) \right\} \right| \\
    & \leq C_n \epsilon^{n+1}\psi_n(x,t;\lambda_i)
  \end{align}
  for $x$ with $-K\leq (-1)^{i-1}x$. The constants $A_{i;n}$ and $B_{i;n}$ are determined from $(M_1,M_2)$ defined by~\eqref{eq:Mi}.
\end{prop}

\begin{rem}
  \label{rem:fn}
  The function $f_{i;n}$ appears in~\cite[Section~4]{vBPW08}. It is shown that $f_{i;n}(z)$ decays exponentially as $(-1)^{i-1}z\to \infty$ but decays algebraically as $f_{i;n}(z)\sim z^{-\alpha_{n-1}}$ in the limit $(-1)^{i-1}z\to -\infty$. In particular, if $|x|\leq K$ for some fixed $K>0$, we have
  \begin{equation*}
    (t+1)^{-\alpha_{n-1}/2}f_{i;n}\left( \frac{x-\lambda_i(t+1)}{\sqrt{t+1}} \right) \sim t^{-\alpha_{n-1}}
  \end{equation*}
\end{rem}

\begin{rem}
  With some additional effort, we can show that, for $-K\leq (-1)^{i-1}x$, the higher-order diffusion waves $\xi_{i;n}$ ($n\geq 1$) are asymptotically equivalent to the higher-order terms of the asymptotic expansion constructed in~\cite{vBPW08}.
\end{rem}

\subsection{Time-asymptotic expansion with pointwise remainder estimates}
Let $\bm{u}_0=(v_0-1,u_0)$ and denote its anti-derivatives by $\bm{u}_{0}^{\pm}$, that is,
\begin{equation}
  \bm{u}_{0}^{-}(x)=\int_{-\infty}^{x}\bm{u}_0(y)\, dy, \quad \bm{u}_{0}^{+}(x)=\int_{x}^{\infty}\bm{u}_0(y)\, dy.
\end{equation}
Our main theorem is the following.

\begin{thm}
  \label{thm:main1}
  For $\bm{u}_0=(v_0-1,u_0)\in H^6(\mathbb{R})\times H^6(\mathbb{R})$, let $(v,u)$ be the solution to~\eqref{eq:p_system}. Define $u_i$, $\theta_i$, and $\xi_{i;n}$ by~\eqref{eq:ui},~\eqref{eq:thetai}, and~\eqref{eq:xii_n}, respectively. Set
  \begin{equation}
    u_{i;1}=\xi_{i;1}+\gamma_{i'}\partial_x \theta_{i'}, \quad u_{i;n}=\xi_{i;n}+\gamma_{i'}\partial_x \xi_{i';n-1} \quad (n\geq 2),
  \end{equation}
  where $i'=3-i$ and $\gamma_i=(-1)^i \nu/(4c)$. Then for $n\geq 1$, there exist positive constants $\delta_n$ and $C_n$ such that if
  \begin{align}
    \label{thm:main1:eq:smallness}
    \begin{aligned}
      \delta
      & \coloneqq \| \bm{u}_0 \|_6+\sup_{x\in \mathbb{R}}[(|x|+1)^{\alpha_n}|\bm{u}_0(x)|+(|x|+1)^{5/4}|\bm{u}_{0}'(x)|] \\
      & \quad +\sup_{x>0}[(|x|+1)^{\beta_n}(|\bm{u}_{0}^{-}(-x)|+|\bm{u}_{0}^{+}(x)|)]\leq \delta_n,
    \end{aligned}
  \end{align}
  the solution $(v,u)$ satisfies the following pointwise estimates:
  \begin{equation}
    \label{thm:main1:eq:PWE}
    \biggl| \biggl( u_i-\theta_i-\sum_{k=1}^{n}u_{i;k} \biggr)(x,t) \biggr| \leq C_n\delta \Psi_{i;n}(x,t)
  \end{equation}
  for all $x\in \mathbb{R}$ and $t\geq 0$. Here $\Psi_{i;n}$ is defined by~\eqref{eq:aux}.
\end{thm}

As a corollary, we obtain the following $L^p$-decay estimates. Combining this with Proposition~\ref{prop:xii_n}, it follows that $u_i \sim \theta_i+\sum_{n=1}^{\infty}\xi_{i;n}$ is a time-asymptotic expansion in the $L^p(\mathbb{R})$-norm for all $1\leq p\leq \infty$.

\begin{cor}
  \label{thm:main1:cor1}
  Under the assumptions of Theorem~\ref{thm:main1}, we have the optimal $L^p$-decay estimate
  \begin{equation}
    \biggl\| \biggl( u_i-\theta_i-\sum_{k=1}^{n}\xi_{i;k} \biggr)(\cdot,t) \biggr\|_{L^p}\leq C_n \delta(t+1)^{-(\alpha_n-1/p)/2} \quad (1\leq p\leq \infty).
  \end{equation}
\end{cor}

\begin{proof}
  The same bound for $u_i-\theta_i-\sum_{k=1}^{n}u_{i;k}$ easily follows from Theorem~\ref{thm:main1}. We can replace $u_{i;k}$ by $\xi_{i;k}$ thanks to~\eqref{eq:theta_explicit} and Proposition~\ref{prop:xii_n}.
\end{proof}

We also obtain the following local-in-space decay estimates.

\begin{cor}
  \label{thm:main1:cor2}
  Under the assumptions of Theorem~\ref{thm:main1}, when $|x|\leq K$ for some fixed $K>0$, we have
  \begin{equation}
    \biggl| \biggl( u_i-\sum_{k=1}^{n}\xi_{i;k} \biggr)(x,t) \biggr|\leq C_{n,K}\delta(t+1)^{-\alpha_n}.
  \end{equation}
  Moreover, there exist constants $\{ A_{i;k} \}_{k=1}^{n}$ determined from $(M_1,M_2)$ such that
  \begin{equation}
    \biggl| u_i(x,t)-\sum_{k=1}^{n}\frac{A_{i;k}}{(t+1)^{\alpha_{k-1}/2}}f_{i;k}\left( \frac{x-\lambda_i(t+1)}{\sqrt{t+1}} \right) \biggr| \leq C_{n,K}\delta(t+1)^{-\alpha_n}.
  \end{equation}
  Here $M_i$ and $f_{i;k}$ are defined by~\eqref{eq:Mi} and~\eqref{eq:fn}, respectively.
\end{cor}

\begin{proof}
  Again, the same bound for $u_i-\theta_i-\sum_{k=1}^{n}u_{i;k}$ easily follows from Theorem~\ref{thm:main1}. We can then replace $u_{i;k}$ by $\xi_{i;k}$ thanks to~\eqref{eq:theta_explicit} and Proposition~\ref{prop:xii_n}. The second inequality follows from Proposition~\ref{prop:xii_n:fn}.
\end{proof}

By Corollary~\ref{thm:main1:cor2}, and also Remark~\ref{rem:fn}, we now have a detailed picture of the power-law asymptotics of the solution around $x=0$ where the diffusion waves decay exponentially.

\begin{rem}
  The term $\gamma_{i'}\partial_x \theta_i$ and $\gamma_{i'}\partial_x \xi_{i;n}$ are both neglected in the two corollaries above. These are negligible in the $L^p(\mathbb{R})$-norm and locally around $x=0$ but are important in the neighborhood of the other characteristic line $x=-\lambda_i t$. For this reason, these terms are required in the statement of Theorem~\ref{thm:main1}.
\end{rem}

\begin{rem}
  The rather strong $H^6$-regularity is required to invoke pointwise estimates of $\partial_x(u_i-\theta_i)$ provided by~\cite[Theorem~2.6 and Remark~2.8]{LZ97}. The proof involves energy estimates up to the $H^6(\mathbb{R})$-norm. These also imply a unique global-in-time existence theorem in appropriate Sobolev spaces. Off course, global-in-time existence of solutions can be proved with much lower regularity~\cite{Kanel68,LY22}, but proving detailed pointwise estimates for such data seems to be difficult at this point.
\end{rem}

\section{Proof}
\label{sec:proof}
The following function appears frequently in the subsequent part of the paper:
\begin{equation}
  \label{eq:Theta}
  \Theta_{\alpha}(x,t;\lambda,\mu)=(t+1)^{-\alpha/2}e^{-\frac{(x-\lambda(t+1))^2}{\mu(t+1)}}.
\end{equation}
Here $\lambda \in \mathbb{R}$ and $\alpha,\mu>0$. Note that
\begin{equation}
  \label{eq:theta_Theta}
  |\theta_i(x,t)|\leq A_0|M_i|\Theta_1(x,t;\lambda_i,2\nu), \quad \Theta_{\alpha_n}(x,t;\lambda,\mu)\leq B_0 \psi_n(x,t;\lambda)
\end{equation}
for some positive constants $A_0$ and $B_0$. In what follows, the symbols $C$ and $\nu^*$ denote sufficiently large constants. 

\subsection{Pointwise estimates of the higher-order diffusion waves}
We start with the proofs of Propositions~\ref{prop:xii_n} and~\ref{prop:xii_n:fn}.

Proposition~\ref{prop:xii_n} follows from the following finer version.

\begin{lem}
  \label{lem:xii_n}
  Let $n\geq 1$ and $\epsilon=\max(M_1,M_2)$. If $\epsilon$ is sufficiently small, we have
  \begin{align}
    \label{lem:xii_n_ineq1}
    \begin{aligned}
      & |\partial_{x}^{k}\xi_{i;n}(x,t)-(-1)^i (2c)^{-1}\partial_{x}^{k}(\theta_{i'}\xi_{i';n-1})(x,t)| \\
      & \quad \leq C_{n,k}\epsilon^{n+1}(t+1)^{-k/2}\psi_{n-1}(x,t;\lambda_i)
    \end{aligned}
  \end{align}
  for any integer $k\geq 0$. In particular, we have
  \begin{align}
    \label{lem:xii_n_ineq2}
    |\partial_{x}^{k}\xi_{i;n}(x,t)|
    & \leq C_{n,k}\epsilon^{n+1}(t+1)^{-k/2}[\psi_{n-1}(x,t;\lambda_i)+\Theta_{2\beta_{n-1}}(x,t;\lambda_{i'},\nu^*)] \\
    & \leq C_{n,k}\epsilon^{n+1}(t+1)^{-k/2}\Psi_{i;n-1}(x,t).
  \end{align}
\end{lem}

\begin{proof}
  We assume $t\geq 4$ in the following (the lemma is otherwise easier to prove). The lemma is trivial for $n=0$ if we set $\xi_{i;0}=\theta_i/2$ and $\xi_{i';-1}=0$. So it suffices to prove the lemma for $n$ assuming that it holds for $n-1\geq 0$. In what follows, we only prove the case of $i=1$ and $k=0$ since the other cases are similar. Note first that by~\eqref{eq:xii_n} and Duhamel's principle, we have $\xi_{i;n}(x,t)=\zeta_{1;n}(x,t)+\eta_{1;n}(x,t)$, where
  \begin{equation}
    \label{eq:zetai_n}
    \zeta_{1;n}(x,t)=
    -(2\pi \nu)^{-1/2}\int_{0}^{t}\int_{-\infty}^{\infty}(t-s)^{-1/2}e^{-\frac{(x-y-c(t-s))^2}{2\nu(t-s)}}\partial_x (\theta_2 \xi_{2;n-1})(y,s)\, dyds
  \end{equation}
  and
  \begin{equation}
    \label{eq:etai_n}
    \eta_{1;n}(x,t)=
    -(2\pi \nu)^{-1/2}\int_{0}^{t}\int_{-\infty}^{\infty}(t-s)^{-1/2}e^{-\frac{(x-y-c(t-s))^2}{2\nu(t-s)}}\partial_x (\theta_1 \xi_{1;n})(y,s)\, dyds.
  \end{equation}

  We first consider $\zeta_{1;n}(x,t)$. Set $I(x,t)=-\sqrt{2\pi \nu}\zeta_{1;n}(x,t)$ and $f=\theta_2 \xi_{2;n-1}$. By Lemma~\ref{lem:heat_equation_technique}, we have
  \begin{equation}
    I(x,t)=(2c)^{-1}\sqrt{2\pi \nu}f(x,t)+I_1(x,t)+I_2(x,t) ,
  \end{equation}
  where
  \begin{equation}
    I_1(x,t)=\int_{0}^{t^{1/2}}\int_{-\infty}^{\infty}\partial_x \left\{ (t-s)^{-1/2}e^{-\frac{(x-y-c(t-s))^2}{2\nu(t-s)}} \right\} f(y,s)\, dyds
  \end{equation}
  and
  \begin{align}
    I_2(x,t)
    & =-(2c)^{-1}\int_{-\infty}^{\infty}(t-t^{1/2})^{-1/2}e^{-\frac{(x-y-c(t-\sqrt{t}))^2}{2\nu(t-\sqrt{t})}}f(y,t^{1/2})\, dy \\
    & \quad -(2c)^{-1}\int_{t^{1/2}}^{t}\int_{-\infty}^{\infty}(t-s)^{-1/2}e^{-\frac{(x-y-c(t-s))^2}{2\nu(t-s)}}L_2 f(y,s)\, dyds \\
    & \eqqcolon I_{21}(x,t)+I_{22}(x,t).
  \end{align}
  Here $L_2=\partial_t-c\partial_x-(\nu/2)\partial_{x}^{2}$. By the induction hypothesis, we have
  \begin{equation}
    |f(x,t)|\leq C\epsilon^{n+1}\Theta_{\alpha_{n-2}+1}(x,t;-c,\nu^*).
  \end{equation}
  By Lemmas~\ref{lem:LZ97_I1} and~\ref{lem:LZ97_I21}, we obtain
  \begin{equation}
    |I_1(x,t)|+|I_{21}(x,t)|\leq C\epsilon^{n+1}\Theta_{\alpha_{n-1}}(x,t;c,\nu^*)\leq C\epsilon^{n+1}\psi_{n-1}(x,t;c).
  \end{equation}
  Next, note that~\eqref{eq:thetai} and~\eqref{eq:xii_n} imply
  \begin{equation}
    |L_2 f(x,t)|\leq C\epsilon^{n+1}\Theta_{\alpha_{n-2}+3}(x,t;-c,\nu^*).
  \end{equation}
  Then by Lemma~\ref{lem:LZ97_I22}, we obtain
  \begin{equation}
    |I_{22}(x,t)|\leq C\epsilon^{n+1}\psi_{n-1}(x,t;c).
  \end{equation}
  We have thus proved
  \begin{equation}
    \label{lem:xii_n:proof:eq1}
    |\zeta_{1;n}(x,t)+(2c)^{-1}(\theta_2 \xi_{2;n-1})(x,t)|\leq C\epsilon^{n+1}\psi_{n-1}(x,t;c).
  \end{equation}

  We next consider $\eta_{1;n}(x,t)$. Note that it is the solution to
  \begin{equation}
    \begin{dcases}
      \partial_t \eta_{1;n}+c\partial_x \eta_{1;n}+\partial_x(\theta_1 \eta_{1;n})=\frac{\nu}{2}\partial_{x}^{2}\eta_{1;n}-\partial_x(\theta_1 \zeta_{1;n}), & x\in \mathbb{R},t>0, \\
      \eta_{1;n}(x,0)=0, & x\in \mathbb{R}.
    \end{dcases}
  \end{equation}
  This variable coefficient equation can be solved by an iteration scheme. Let $\eta_{1;n}^{(1)}$ be the solution to
  \begin{equation}
    \begin{dcases}
      \partial_t \eta_{1;n}^{(1)}+c\partial_x \eta_{1;n}^{(1)}=\frac{\nu}{2}\partial_{x}^{2}\eta_{1;n}^{(1)}-\partial_x (\theta_1 \zeta_{1;n}), & x\in \mathbb{R},t>0, \\
      \eta_{1;n}^{(1)}(x,0)=0, & x\in \mathbb{R}
    \end{dcases}
  \end{equation}
  and $\eta_{1;n}^{(k)}$ ($k\geq 2$) be the solution to
  \begin{equation}
    \begin{dcases}
      \partial_t \eta_{1;n}^{(k)}+c\partial_x \eta_{1;n}^{(k)}=\frac{\nu}{2}\partial_{x}^{2}\eta_{1;n}^{(k)}-\partial_x(\theta_1 \eta_{1;n}^{(k-1)}) & x\in \mathbb{R},t>0, \\
      \eta_{1;n}^{(k)}(x,0)=0, & x\in \mathbb{R}.
    \end{dcases}
  \end{equation}
  Then we can write $\eta_{1;n}$ as
  \begin{equation}
    \label{lem:xii_n:proof:eq11}
    \eta_{1;n}(x,t)=\sum_{k=1}^{\infty}\eta_{1;n}^{(k)}(x,t).
  \end{equation}
  We now give bounds for $\eta_{1;n}^{(k)}$ ($k\geq 1$) inductively. Note first that
  \begin{equation}
    \label{eq:eta_integral_rep}
    \eta_{1;n}^{(1)}(x,t)=-(2\pi \nu)^{-1/2}\int_{0}^{t}\int_{-\infty}^{\infty}(t-s)^{-1/2}e^{-\frac{(x-y-c(t-s))^2}{2\nu(t-s)}}\partial_x (\theta_1 \zeta_{1;n})(y,s)\, dyds
  \end{equation}
  and that~\eqref{lem:xii_n:proof:eq1} implies
  \begin{equation}
    |(\theta_1 \zeta_{1;n})(x,t)|\leq A_1 \epsilon^{n+2}\Theta_{\alpha_{n-1}+1}(x,t;c,\nu')
  \end{equation}
  for some positive constants $A_1$ and $\nu'$. Then by~\cite[Lemma~3.2]{LZ97}, we obtain
  \begin{equation}
    |\eta_{1;n}^{(1)}(x,t)|\leq M A_1 \epsilon^{n+2}\Theta_{\alpha_{n-1}}(x,t;c,\nu')
  \end{equation}
  for some $M>0$. This means that the inequality
  \begin{equation}
    \label{lem:xii_n:proof:eq2}
    |\eta_{1;n}^{(l)}(x,t)|\leq MA_1 (MA_0)^{l-1}\epsilon^{n+l+1}\Theta_{\alpha_{n-1}}(x,t;c,\nu')
  \end{equation}
  holds for $l=1$. We then show that~\eqref{lem:xii_n:proof:eq2} holds for $l=k+1$ assuming that it holds for $l=k$. By the induction hypothesis and~\eqref{eq:theta_Theta}, we have
  \begin{equation}
    |(\theta_1 \eta_{1;n}^{(k)})(x,t)|\leq A_1 (MA_0)^{k}\epsilon^{n+k+2}\Theta_{\alpha_{n-1}+1}(x,t;c,\nu').
  \end{equation}
  Applying~\cite[Lemma~3.2]{LZ97} again, this time to the integral representation
  \begin{equation}
    \eta_{1;n}^{(k+1)}(x,t)=-(2\pi \nu)^{-1/2}\int_{0}^{t}\int_{-\infty}^{\infty}(t-s)^{-1/2}e^{-\frac{(x-y-c(t-s))^2}{2\nu(t-s)}}\partial_x (\theta_1 \eta_{1;n}^{(k)})(y,s)\, dyds,
  \end{equation}
  we obtain
  \begin{equation}
    |\eta_{1;n}^{(k+1)}(x,t)|\leq MA_1 (MA_0)^k \epsilon^{n+k+2}\Theta_{\alpha_{n-1}}(x,t;c,\nu').
  \end{equation}
  Therefore,~\eqref{lem:xii_n:proof:eq2} holds for any $l\geq 1$, and by taking $\epsilon$ sufficiently small, we get
  \begin{equation}
    \label{eq:etai_n:eq1}
    |\eta_{1;n}(x,t)|=\biggl| \sum_{k=1}^{\infty}\eta_{1;n}^{(k)}(x,t) \biggr|\leq C\epsilon^{n+2}\Theta_{\alpha_{n-1}}(x,t;c,\nu')\leq C\epsilon^{n+2}\psi_{n-1}(x,t;c).
  \end{equation}
  Combining this with~\eqref{lem:xii_n:proof:eq1}, we obtain~\eqref{lem:xii_n_ineq1}.
\end{proof}

\begin{rem}
  \label{rem:uniform_in_m_geq_n}
  The proof above can be modified to show that
  \begin{align}
    \begin{aligned}
      & |\partial_{x}^{k}\xi_{i;m}(x,t)-(-1)^i (2c)^{-1}\partial_{x}^{k}(\theta_{i'}\xi_{i';m-1})(x,t)| \\
      & \quad \leq C_{n,k}\epsilon^{m+1}(t+1)^{-k/2}\psi_{n-1}(x,t;\lambda_i)
    \end{aligned}
  \end{align}
  holds for all $m\geq n$ with the smallness of $\epsilon$ depending only on $n$ and $k$.
\end{rem}

We next prove Proposition~\ref{prop:xii_n:fn}. (The proof is rather lengthy and may be skipped; the rest of the paper can be read independently.) Define $\zeta_{i;n}$ and $\eta_{i;n}$ by~\eqref{eq:zetai_n} and~\eqref{eq:etai_n}, respectively. Then Proposition~\ref{prop:xii_n:fn} is a direct consequence of the following lemma.

\begin{lem}
  \label{lem:xi_n:etai_n}
  Let $n\geq 1$ and $\epsilon=\max(M_1,M_2)$. Fix $k\geq 0$. For any $K>0$, if $\epsilon$ is sufficiently small, there exist $A_{i;n}$, $B_{i;n}$, and $C_{n,k}>0$ such that
  \begin{equation}
    \label{lem:xi_n:etai_n:eq:zetai_n}
    \left| \partial_{x}^{k}\left\{ \zeta_{i;n}(x,t)-\frac{A_{i;n}}{(t+1)^{\alpha_{n-1}/2}}f_{i;n}\left( \frac{x-\lambda_i(t+1)}{\sqrt{t+1}} \right) \right\} \right| \leq C_{n,k}\epsilon^{n+1}(t+1)^{-k/2}\psi_n(x,t;\lambda_i)
  \end{equation}
  and
  \begin{equation}
    \label{lem:xi_n:etai_n:eq:etai_n}
    \left| \partial_{x}^{k}\left\{ \eta_{i;n}(x,t)-\frac{B_{i;n}}{(t+1)^{\alpha_{n-1}/2}}g\left( \frac{x-\lambda_i(t+1)}{\sqrt{t+1}} \right) \right\} \right| \leq C_{n,k}\epsilon^{n+2}(t+1)^{-k/2}\psi_n(x,t;\lambda_i)
  \end{equation}
  when $-K\leq (-1)^{i-1}x$. Here $g$ and $f_{i;n}$ are defined by~\eqref{eq:g} and~\eqref{eq:fn}, respectively.
\end{lem}

\begin{proof}
  The lemma is proved by induction in $n$.

  We first consider the case of $n=1$. Let $i=1$ and $k=0$ (the other cases are similar). We start with the proof of~\eqref{lem:xi_n:etai_n:eq:zetai_n}. Note that~\eqref{eq:zetai_n} implies
  \begin{equation}
    \zeta_{1;1}(x,t)=-\frac{1}{2\sqrt{2\pi \nu}}\partial_x \int_{0}^{t}\int_{-\infty}^{\infty}(t-s)^{-1/2}e^{-\frac{(x-y-c(t-s))^2}{2\nu(t-s)}}\theta_{2}^{2}(y,s)\, dyds.
  \end{equation}
  In addition, by~\eqref{eq:theta_explicit} and~\eqref{eq:f0}, we have
  \begin{equation}
    \theta_i(x,t)=(t+1)^{-1/2}f_{i;0}\left( \frac{x-\lambda_i (t+1)}{\sqrt{t+1}} \right).
  \end{equation}
  Hence we may write
  \begin{equation}
    \theta_{2}^{2}(x,t)=\frac{a_{1;1}}{(t+1)\sqrt{2\pi \nu}}e^{-\frac{(x+c(t+1))^2}{2\nu(t+1)}}+\partial_x r(x,t),
  \end{equation}
  where
  \begin{equation}
    a_{1;1}=\int_{-\infty}^{\infty}f_{2;0}^{2}(z)\, dz
  \end{equation}
  and
  \begin{equation}
    r(x,t)=\int_{-\infty}^{x}\left( \theta_{2}^{2}(z,t)-\frac{a_{1;1}}{(t+1)\sqrt{2\pi \nu}}e^{-\frac{(z+c(t+1))^2}{2\nu(t+1)}} \right) \, dz.
  \end{equation}
  Noting that $\lim_{x\to \infty}r(x,t)=0$, we can show that
  \begin{equation}
    |r(x,t)|\leq C\epsilon^2 \Theta_2(x,t;-c,\nu^*), \quad |L_2 r(x,t)|\leq C\epsilon^2 \Theta_4(x,t;-c,\nu^*),
  \end{equation}
  where $L_2=\partial_t-c\partial_x-(\nu/2)\partial_{x}^{2}$. We then have
  \begin{align}
    \zeta_{1;1}(x,t)
    & =-\frac{a_{1;1}}{4\pi \nu}\partial_x \int_{0}^{t}\int_{-\infty}^{\infty}(t-s)^{-1/2}(s+1)^{-1}e^{-\frac{(x-y-c(t-s))^2}{2\nu(t-s)}}e^{-\frac{(y+c(s+1))^2}{2\nu(s+1)}}\, dyds \\
    & \quad -\frac{1}{2\sqrt{2\pi \nu}}\partial_x \int_{0}^{t}\int_{-\infty}^{\infty}(t-s)^{-1/2}e^{-\frac{(x-y-c(t-s))^2}{2\nu(t-s)}}\partial_y r(y,s)\, dyds \\
    & =-\frac{a_{1;1}}{2\sqrt{2\pi \nu}}\partial_x \int_{0}^{t}(t+1)^{-1/2}(s+1)^{-1/2}e^{-\frac{(x-c(t-s)+c(s+1))^2}{2\nu(t+1)}}\, ds \\
    & \quad -\frac{1}{2\sqrt{2\pi \nu}}\partial_x \int_{0}^{t}\int_{-\infty}^{\infty}(t-s)^{-1/2}e^{-\frac{(x-y-c(t-s))^2}{2\nu(t-s)}}\partial_y r(y,s)\, dyds.
  \end{align}
  Concerning the second term, similar calculations leading to the bound of $\zeta_{1;n}(x,t)$ in Lemma~\ref{lem:xii_n} imply
  \begin{equation}
    \left| \partial_x \int_{0}^{t}\int_{-\infty}^{\infty}(t-s)^{-1/2}e^{-\frac{(x-y-c(t-s))^2}{2\nu(t-s)}}\partial_y r(y,s)\, dyds \right| \leq C\epsilon^2 \psi_n(x,t;c)
  \end{equation}
  for $-K\leq x$. For the first term, note that a simple change of variable yields
  \begin{equation}
    -\frac{(t+1)^{-3/4}}{\nu \sqrt{2c}}f_{1;1}\left( \frac{x-c(t+1)}{\sqrt{t+1}} \right)=\partial_x \int_{-1}^{\infty}(t+1)^{-1/2}(s+1)^{-1/2}e^{-\frac{(x-c(t-s)+c(s+1))^2}{2\nu(t+1)}}\, ds.
  \end{equation}
  Therefore,
  \begin{align}
    & \zeta_{1;1}(x,t)-\frac{a_{1;1}}{2\nu \sqrt{4\pi \nu c}}(t+1)^{-3/4}f_{1;1}\left( \frac{x-c(t+1)}{\sqrt{t+1}} \right) \\
    & =O(\epsilon^2)\partial_x \int_{s\in (-1,0)\cup (t,\infty)}(t+1)^{-1/2}(s+1)^{-1/2}e^{-\frac{(x-c(t-s)+c(s+1))^2}{2\nu(t+1)}}\, ds \\
    & \quad +O(\epsilon^2)\psi_n(x,t;c).
  \end{align}
  We then set
  \begin{equation}
    I(x,t)=\partial_x \int_{s\in (-1,0)\cup (t,\infty)}(t+1)^{-1/2}(s+1)^{-1/2}e^{-\frac{(x-c(t-s)+c(s+1))^2}{2\nu(t+1)}}\, ds
  \end{equation}
  and show that $|I(x,t)|\leq C\Theta_2(x,t;c,\nu^*)$ for $-K\leq x$. We first consider Case (i) $|x-c(t+1)|\leq (t+1)^{1/2}$. In this case, we simply have
  \begin{equation}
    |I(x,t)|\leq C(t+1)^{-1}\leq C\Theta_2(x,t;c,\nu^*).
  \end{equation}
  We next consider Case (ii) $-K\leq x\leq c(t+1)-(t+1)^{1/2}$. The integral over $(-1,0)$ is easy to handle. For $s\in (t,\infty)$ on the other hand, when $t$ is large (the case when $t$ is not large is easier), we have
  \begin{gather}
    0\leq c(t+1)-x-2K\leq x-c(t-s)+c(s+1), \\
    0\leq c(s+1)-K\leq x-c(t-s)+c(s+1).
  \end{gather}
  Hence
  \begin{align}
    & \left| \partial_x \int_{t}^{\infty}(t+1)^{-1/2}(s+1)^{-1/2}e^{-\frac{(x-c(t-s)+c(s+1))^2}{2\nu(t+1)}}\, ds \right| \\
    & \quad \leq C\int_{0}^{\infty}(s+1)^{-1/2}e^{-\frac{s^2}{C(t+1)}}\, ds\cdot \Theta_2(x,t;c,\nu^*)\leq C\Theta_2(x,t;c,\nu^*).
  \end{align}
  We end the analysis of $\zeta_{1;1}$ by considering Case (iii) $x\geq c(t+1)+(t+1)^{1/2}$. When $s>-1$, we have
  \begin{gather}
    0\leq x-c(t+1)\leq x-c(t-s)+c(s+1)=x-c(t+1)+2c(s+1), \\
    0\leq 2c(s+1)\leq x-c(t-s)+c(s+1).
  \end{gather}
  From these, it follows that $|I(x,t)|\leq C\Theta_2(x,t;c,\nu^*)$ as in Case (ii). These prove~\eqref{lem:xi_n:etai_n:eq:zetai_n} for $n=1$ by setting
  \begin{equation}
    A_{1;1}=\frac{a_{1;1}}{2\nu \sqrt{4\pi \nu c}}.
  \end{equation}
  We next prove~\eqref{lem:xi_n:etai_n:eq:etai_n} by using the series representation~\eqref{lem:xii_n:proof:eq11}. We first consider
  \begin{equation}
    \eta_{1;1}^{(1)}(x,t)=-(2\pi \nu)^{-1/2}\int_{0}^{t}\int_{-\infty}^{\infty}(t-s)^{-1/2}e^{-\frac{(x-y-c(t-s))^2}{2\nu(t-s)}}\partial_x (\theta_1 \zeta_{1;1})(y,s)\, dyds.
  \end{equation}
  The bound~\eqref{lem:xi_n:etai_n:eq:zetai_n} for $\zeta_{1;1}(x,t)$ implies
  \begin{equation}
    (\theta_1 \zeta_{1;1})(x,t)=\theta_1(x,t)\frac{A_{1;1}}{(t+1)^{3/4}}f_{1;1}\left( \frac{x-c(t+1)}{\sqrt{t+1}} \right)+O(\epsilon^3)\Theta_{\alpha_n+1}(x,t;c,\nu^*).
  \end{equation}
  This holds for all $x\in \mathbb{R}$ since $\theta_1(x,t)$ decays exponentially for $x\leq -K$. Plugging this into~\eqref{eq:eta_integral_rep} and arguing similarly to the analysis of $\zeta_{1;1}$, we get
  \begin{align}
    \eta_{1;1}^{(1)}(x,t)
    & =-\frac{A_{1;1}b_{1;1}}{2\pi \nu}\partial_x \int_{0}^{t}\int_{-\infty}^{\infty}(t-s)^{-1/2}(s+1)^{-5/4}e^{-\frac{(x-y-c(t-s))^2}{2\nu(t-s)}}e^{-\frac{(y-c(s+1))^2}{2\nu(s+1)}}\, dyds \\
    & \quad +O(\epsilon^3)\psi_n(x,t;c) \\
    & =-\frac{A_{1;1}b_{1;1}}{\sqrt{2\pi \nu}}\partial_x \int_{0}^{t}(t+1)^{-1/2}(s+1)^{-3/4}e^{-\frac{(x-c(t+1))^2}{2\nu(t+1)}}\, ds+O(\epsilon^3)\psi_n(x,t;c) \\
    & =-\frac{4A_{1;1}b_{1;1}}{\sqrt{2\pi \nu}}(t+1)^{-3/4}g\left( \frac{x-c(t+1)}{\sqrt{t+1}} \right)+O(\epsilon^3)\psi_n(x,t;c),
  \end{align}
  where
  \begin{equation}
    b_{1;1}=\int_{-\infty}^{\infty}(f_{1;0}f_1)(z)\, dz.
  \end{equation}
  Similar analysis for $\eta_{1;1}^{(k)}(x,t)$ ($k\geq 2$) shows that
  \begin{equation}
    \eta_{1;1}^{(k)}(x,t)=O(\epsilon^{k+2})(t+1)^{-3/4}g\left( \frac{x-c(t+1)}{\sqrt{t+1}} \right)+O(\epsilon^{k+2})\psi_n(x,t;c).
  \end{equation}
  Taking the sum $\sum_{k=1}^{\infty}$, it follows that~\eqref{lem:xi_n:etai_n:eq:etai_n} holds for $n=1$ with
  \begin{equation}
    B_{1;1}=-\frac{4A_{1;1}b_{1;1}}{\sqrt{2\pi \nu}}+O(\epsilon^4).
  \end{equation}

  We next prove the lemma for $n$ assuming that it holds for $n-1$. Let $i=1$ and $k=0$ (the other cases are similar). The induction hypothesis and Lemma~\ref{lem:xii_n} imply
  \begin{align}
    & h(x,t)\coloneqq (\theta_2 \xi_{2;n-1})(x,t)-\frac{A_{2;n-1}}{(t+1)^{\alpha_{n-2}/2}}\theta_2(x,t)f_{2;n-1}\left( \frac{x+c(t+1)}{\sqrt{t+1}} \right) \\
    & \quad -\frac{B_{2;n-1}}{(t+1)^{\alpha_{n-2}/2}}\theta_2(x,t)g\left( \frac{x+c(t+1)}{\sqrt{t+1}} \right)=O(\epsilon^{n+1})\Theta_{\alpha_{n-1}+1}(x,t;-c,\nu^*)
  \end{align}
  for all $x\in \mathbb{R}$ (not just for $x\leq K$). We also have $\partial_x h(x,t)=O(\epsilon^{n+1})\Theta_{\alpha_{n-1}+2}(x,t;-c,\nu^*)$. Using these, we can show that
  \begin{equation}
    L_2 h(x,t)=O(\epsilon^{n+1})\Theta_{\alpha_{n-1}+3}(x,t;-c,\nu^*),
  \end{equation}
  where $L_2=\partial_t-c\partial_x-(\nu/2)\partial_{x}^{2}$. Then similar calculations leading to the bound of $\zeta_{1;1}(x,t)$ above imply~\eqref{lem:xi_n:etai_n:eq:zetai_n} with
  \begin{equation}
    A_{1;n}=\frac{1}{\nu \sqrt{4\pi \nu c}}(A_{2;n-1}a_{1;n}+B_{2;n-1}b_{1;n}),
  \end{equation}
  where
  \begin{equation}
    a_{1;n}=\int_{-\infty}^{\infty}(f_{2;0}f_{2;n-1})(-z)\, dz, \quad b_{1;n}=\int_{-\infty}^{\infty}(gf_{2;0})(z)\, dz.
  \end{equation}
  The bound~\eqref{lem:xi_n:etai_n:eq:etai_n} for $\eta_{1;n}(x,t)$ is proved in a way similar to that for $n=1$. This ends the proof of the lemma.
\end{proof}

For the proof of Theorem~\ref{thm:main1}, it is convenient to unify $(\xi_{i;n})_{n=1}^{\infty}$ into a single function
\begin{equation}
  \label{eq:Xii}
  \Xi_i(x,t)=\sum_{n=1}^{\infty}\xi_{i;n}(x,t).
\end{equation}
Taking the infinite sum of~\eqref{eq:xii_n}, we see that $(\Xi_1,\Xi_2)$ is the solution to the system
\begin{equation}
  \label{eq:Xi}
  \begin{dcases}
    \partial_t \Xi_1+c\partial_x \Xi_1+\partial_x(\theta_{2}^{2}/2+\theta_1 \Xi_1+\theta_2 \Xi_2)=\frac{\nu}{2}\partial_{x}^{2}\Xi_1, & x\in \mathbb{R},t>0, \\
    \partial_t \Xi_2-c\partial_x \Xi_2+\partial_x(\theta_{1}^{2}/2+\theta_1 \Xi_1+\theta_2 \Xi_2)=\frac{\nu}{2}\partial_{x}^{2}\Xi_2, & x\in \mathbb{R},t>0, \\
    \Xi_1(x,0)=\Xi_2(x,0)=0, & x\in \mathbb{R}.
  \end{dcases}
\end{equation}
Then Lemma~\ref{lem:xii_n} and Remark~\ref{rem:uniform_in_m_geq_n} imply the following.

\begin{lem}
  \label{lem:Xii}
  Let
  \begin{equation}
    \Xi_{i;n}(x,t)=\sum_{m=n+1}^{\infty}\xi_{i;m}(x,t) \quad (n\geq -1).
  \end{equation}
  Here $\xi_{i;0}=\theta_i/2$. Then for $n\geq 0$, if $\epsilon=\max(M_1,M_2)$ is sufficiently small, we have
  \begin{equation}
    |\partial_{x}^{k}\Xi_{i;n}(x,t)-(-1)^i (2c)^{-1}\partial_{x}^{k}(\theta_{i'}\Xi_{i';n-1})(x,t)|\leq C_{n,k}\epsilon^{n+2}(t+1)^{-k/2}\psi_n(x,t;\lambda_i)
  \end{equation}
  for any integer $k\geq 0$. In particular, we have
  \begin{equation}
    |\partial_{x}^{k}\Xi_{i;n}(x,t)|\leq C_{n,k}\epsilon^{n+2}(t+1)^{-k/2}\Psi_{i;n}(x,t).
  \end{equation}
\end{lem}

\subsection{Proof of Theorem~\ref{thm:main1}}
Let us explain the strategy to prove Theorem~\ref{thm:main1}. Note first that by Lemma~\ref{lem:Xii}, it suffices to prove the following.

\begin{thm}
  \label{thm:main2}
  For $\bm{u}_0=(v_0-1,u_0)\in H^6(\mathbb{R})\times H^6(\mathbb{R})$, let $(v,u)$ be the solution to~\eqref{eq:p_system}. Define $u_i$, $\theta_i$, and $\Xi_i$ by~\eqref{eq:ui},~\eqref{eq:thetai}, and~\eqref{eq:Xi}, respectively. Then for $n\geq 1$, there exist positive constants $\delta_n$ and $C_n$ such that if $\delta \leq \delta_n$, where $\delta$ is defined by~\eqref{thm:main1:eq:smallness}, the solution $(v,u)$ satisfies the following pointwise estimates:
  \begin{equation}
    \label{thm:main2:eq:PWE}
    |(u_i-\theta_i-\Xi_i-\gamma_{i'}\partial_x \theta_{i'}-\gamma_{i'}\partial_x \Xi_{i'})(x,t)|\leq C_n\delta \Psi_{i;n}(x,t)
  \end{equation}
  for all $x\in \mathbb{R}$ and $t\geq 0$. Here $\gamma_i=(-1)^i \nu/(4c)$ and $i'=3-i$.
\end{thm}

To prove Theorem~\ref{thm:main2}, we set
\begin{equation}
  \label{eq:vi}
  v_i=u_i-\theta_i-\Xi_i-\gamma_{i'}\partial_x \theta_{i'}-\gamma_{i'}\partial_x \Xi_{i'}
\end{equation}
and define $P(t)$ by
\begin{equation}
  \label{eq:P}
  P(t)\coloneqq \sum_{i=1}^{2}\sup_{0\leq s\leq t}\left| v_i(\cdot,s)\Psi_{i;n}(\cdot,s)^{-1} \right|_{L^{\infty}}.
\end{equation}
Our goal is then to prove the inequality
\begin{equation}
  \label{eq:P_ineq}
  P(t)\leq C\delta+C(\delta+P(t))^2 \quad (t\geq 0).
\end{equation}
From this inequality, taking $\delta$ sufficiently small, we can conclude that $P(t)\leq C\delta$ for all $t\geq 0$ by a standard argument (see Section~\ref{sec:final}).

\begin{rem}
  For the argument above to work, we first need to show that $P(t)$ is finite. This can be proved, for example, by examining the iterative scheme in~\cite[Section~2.1]{Koike21} for the construction of the local-in-time solution to~\eqref{eq:p_system}. The key step of the scheme consists of solving a variable coefficient parabolic equation, and by the Levi parametrix method, we can prove a gaussian upper bound for the fundamental solution. This bound allows us to control the spatial decay of each approximate solution, and by taking the limit, we can check that $P(t)$ is finite at least for a short period of time. By the calculations below, it follows that~\eqref{eq:P_ineq} and hence $P(t)\leq C\delta$ hold for this short duration. Then a standard continuity argument shows that $P(t)\leq C\delta$ actually holds for all $t\geq 0$.
\end{rem}

The proof of~\eqref{eq:P_ineq} is based on pointwise estimates of Greens' function~\cite{LZ97,LZ09} which we shall explain in the next section. We also give an integral formulation of~\eqref{eq:p_system}. In the remaining sections, we prove bounds for the terms appearing in the integral equations which yield~\eqref{eq:P_ineq}.

\subsubsection{Pointwise estimates of Green's function and integral equations}
Our equations~\eqref{eq:p_system} can be written in the form
\begin{equation}
  \label{eq:p_system:vector_form}
  \bm{u}_t+A\bm{u}_x=
  \begin{pmatrix}
    0 & 0 \\
    0 & \nu
  \end{pmatrix}
  \bm{u}_{xx}+
  \begin{pmatrix}
    0 \\
    N_x
  \end{pmatrix}
\end{equation}
with
\begin{equation}
  \label{eq:N}
  \bm{u}=
  \begin{pmatrix}
    v-1 \\
    u
  \end{pmatrix},
  \quad A=
  \begin{pmatrix}
    0 & -1 \\
    -c^2 & 0
  \end{pmatrix},
  \quad N=-p(v)+p(1)-c^2(v-1)-\nu \frac{v-1}{v}u_x.
\end{equation}
The matrix $A$ has right and left eigenvectors $r_i$ and $l_i$ ($i=1,2$), corresponding to the eigenvalue $\lambda_i=(-1)^{i-1}c$, given by
\begin{equation}
  r_i=\frac{2c}{p''(1)}
  \begin{bmatrix}
    (-1)^i \\
    c
  \end{bmatrix},
  \quad l_i=\frac{p''(1)}{4c}
  \begin{bmatrix}
    (-1)^i & 1/c
  \end{bmatrix}.
\end{equation}
We note that~\eqref{eq:ui} can be written as $u_i=l_i(v-1 \, u)^{T}$.

We define Green's function $G=G(x,t)\in \mathbb{R}^{2\times 2}$ for the linearization of~\eqref{eq:p_system:vector_form} as the solution to
\begin{equation}
  \begin{dcases}
    \partial_t G+A\partial_x G=
    \begin{pmatrix}
      0 & 0 \\
      0 & \nu
    \end{pmatrix}
    \partial_{x}^{2}G,   & x\in \mathbb{R},\, t>0, \\
    G(x,0)=\delta(x)I_2, & x\in \mathbb{R},
  \end{dcases}
\end{equation}
where $\delta(x)$ is the Dirac delta function and $I_2$ is the $2\times 2$ identity matrix. In addition, define $G^*=G^*(x,t)\in \mathbb{R}^{2\times 2}$ by
\begin{equation}
  G^*(x,t)=\frac{1}{2(2\pi \nu t)^{1/2}}e^{-\frac{(x-ct)^2}{2\nu t}}
  \begin{pmatrix}
    1  & -1/c \\
    -c & 1
  \end{pmatrix}
  +\frac{1}{2(2\pi \nu t)^{1/2}}e^{-\frac{(x+ct)^2}{2\nu t}}
  \begin{pmatrix}
    1 & 1/c \\
    c & 1
  \end{pmatrix}.
\end{equation}

The following theorem is of fundamental importance in our analysis.

\begin{thm}[{\cite[Theorem~5.8]{LZ97} and~\cite[Theorem~1.3]{LZ09}}]
  \label{thm:G_PWE}
  For any $k\geq 0$, we have
  \begin{equation}
    \biggl| \partial_{x}^{k}G(x,t)-\partial_{x}^{k}G^*(x,t)-e^{-\frac{c^2}{\nu}t}\sum_{j=0}^{k}\delta^{(k-j)}(x)Q_j(t) \biggr|\leq C(t+1)^{-\frac{1}{2}}t^{-\frac{k+1}{2}}\sum_{i=1}^{2}e^{-\frac{(x-\lambda_i t)^2}{Ct}},
  \end{equation}
  where $\delta^{(k)}(x)$ is the $k$-th derivative of the Dirac delta function and $Q_j=Q_j(t)$ is a $2\times 2$ polynomial matrix. In particular,
  \begin{equation}
    Q_0=
    \begin{pmatrix}
      1 & 0 \\
      0 & 0
    \end{pmatrix}
    , \quad Q_1=
    \begin{pmatrix}
      0        & -1/\nu \\
      -c^2/\nu & 0
    \end{pmatrix}.
  \end{equation}
  Moreover, with $\gamma_i=(-1)^i \nu/(4c)$, we have
  \begin{align}
    \begin{aligned}
      & \partial_{x}^{k}G(x,t)-\partial_{x}^{k}G^*(x,t)-\partial_{x}^{k+1}\sum_{i=1}^{2}\gamma_i \frac{e^{-\frac{(x-\lambda_i t)^2}{2\nu t}}}{(2\pi \nu t)^{1/2}}
      \begin{pmatrix}
        -1 & 0 \\
        0 & 1
      \end{pmatrix}
      -e^{-\frac{c^2}{\nu}t}\sum_{j=0}^{k}\delta^{(k-j)}(x)Q_j(t) \\
      & =O(1)(t+1)^{-\frac{1}{2}}t^{-\frac{k+1}{2}}e^{-\frac{(x-ct)^2}{Ct}}
      \begin{pmatrix}
        1  & -1/c \\
        -c & 1
      \end{pmatrix}
      +O(1)(t+1)^{-\frac{1}{2}}t^{-\frac{k+1}{2}}e^{-\frac{(x+ct)^2}{Ct}}
      \begin{pmatrix}
        1 & 1/c \\
        c & 1
      \end{pmatrix} \\
      & \quad +O(1)(t+1)^{-\frac{1}{2}}t^{-\frac{k+2}{2}}\sum_{i=1}^{2}e^{-\frac{(x-\lambda_i t)^2}{Ct}}.
    \end{aligned}
  \end{align}
   Here $O(1)$ is a bounded scalar function.
\end{thm}

For the analysis of $u_i$, we need pointwise estimates for
\begin{equation}
  g_i=
  \begin{pmatrix}
    g_{i1} & g_{i2}
  \end{pmatrix}
  \coloneqq l_i G
  \begin{pmatrix}
    r_1 & r_2
  \end{pmatrix},
  \quad g_{i}^{*}=
  \begin{pmatrix}
    g_{i1}^{*} & g_{i2}^{*}
  \end{pmatrix}
  \coloneqq l_i G^*
  \begin{pmatrix}
    r_1 & r_2
  \end{pmatrix}.
\end{equation}
We note that
\begin{equation}
  \label{eq:gij_aster}
  g_{ij}^{*}=(2\pi \nu t)^{-1/2}e^{-\frac{(x-\lambda_i t)^2}{2\nu t}}\delta_{ij},
\end{equation}
where $\delta_{ij}$ is the Kronecker delta. Then Theorem~\ref{thm:G_PWE} implies
\begin{equation}
  \label{eq:gi_PWE1}
  \biggl| \partial_{x}^{k}g_i(x,t)-\partial_{x}^{k}g_{i}^{*}(x,t)-e^{-\frac{c^2}{\nu}t}\sum_{j=0}^{k}\delta^{(k-j)}(x)q_{ik}(t) \biggr| \leq C(t+1)^{-\frac{1}{2}}t^{-\frac{k+1}{2}}\sum_{i=1}^{2}e^{-\frac{(x-\lambda_i t)^2}{Ct}},
\end{equation}
where
\begin{equation}
  q_{ik}(t)=l_i Q_k(t)
  \begin{pmatrix}
    r_1 & r_2
  \end{pmatrix}.
\end{equation}
Moreover, we have
\begin{align}
  \label{eq:gi_PWE2}
  \begin{aligned}
    & \partial_{x}^{k}g_i(x,t)-\partial_{x}^{k}g_{i}^{*}(x,t)-\gamma_{i'}\partial_{x}^{k+1}g_{i'}^{*}(x,t)-e^{-\frac{c^2}{\nu}t}\sum_{j=0}^{k}\delta^{(k-j)}(x)q_{ik}(t) \\
    & \quad =O(1)(t+1)^{-\frac{1}{2}}t^{-\frac{k+1}{2}}e^{-\frac{(x-ct)^2}{Ct}}+O(1)(t+1)^{-\frac{1}{2}}t^{-\frac{k+2}{2}}e^{-\frac{(x+ct)^2}{Ct}}.
  \end{aligned}
\end{align}

We next write down an integral equation for $v_i$ defined by~\eqref{eq:vi}. Let
\begin{equation}
  n=l_i
  \begin{pmatrix}
    0 \\
    N
  \end{pmatrix}
  =\frac{p''(1)}{4c^2}N, \quad n^*=-\theta_{1}^{2}/2-\theta_{2}^{2}/2-\theta_1 \Xi_1-\theta_2 \Xi_2,
\end{equation}
where $\Xi_i$ and $N$ are defined by~\eqref{eq:Xii} and~\eqref{eq:N}, respectively. Then by Duhamel's principle, we obtain the following.

\begin{lem}
  \label{lem:vi_IE}
  The function $v_i$ defined by~\eqref{eq:vi} satisfies the integral equation
  \begin{align}
    \label{prop:vi_IE:eq:vi_IE}
    \begin{aligned}
      v_i(x,t)
      & =\int_{-\infty}^{\infty}g_i(x-y,t)
      \begin{pmatrix}
        u_1 \\
        u_2
      \end{pmatrix}
      (y,0)\, dy-\int_{-\infty}^{\infty}g_{i}^{*}(x-y,t)
      \begin{pmatrix}
        \theta_1 \\
        \theta_2
      \end{pmatrix}
      (y,0)\, dy \\
      & \quad -\gamma_{i'}\int_{-\infty}^{\infty}\partial_x g_{i'}^{*}(x-y,t)
      \begin{pmatrix}
        \theta_1 \\
        \theta_2
      \end{pmatrix}
      (y,0)\, dy \\
      & \quad +\int_{0}^{t}\int_{-\infty}^{\infty}g_{ii}^{*}(x-y,t-s)\partial_x (n-n^*)(y,s)(y,s)\, dyds \\
      & \quad +\sum_{j=1}^{2}\int_{0}^{t}\int_{-\infty}^{\infty}(g_{ij}-g_{ij}^{*})(x-y,t-s)\partial_x n(y,s)\, dyds \\
      & \quad -\gamma_{i'}\int_{0}^{t}\int_{-\infty}^{\infty}\partial_x g_{i'i'}^{*}(x-y,t-s)\partial_x n^*(y,s)(y,s)\, dyds.
    \end{aligned}
  \end{align}
  Here $\gamma_i=(-1)^i \nu/(4c)$ and $i'=3-i$.
\end{lem}

We set
\begin{align}
  \label{eq:Ii}
  \begin{aligned}
    \mathcal{I}_i(x,t)
    & =\int_{-\infty}^{\infty}g_i(x-y,t)
    \begin{pmatrix}
      u_1 \\
      u_2
    \end{pmatrix}
    (y,0)\, dy-\int_{-\infty}^{\infty}g_{i}^{*}(x-y,t)
    \begin{pmatrix}
      \theta_1 \\
      \theta_2
    \end{pmatrix}
    (y,0)\, dy \\
    & \quad -\gamma_{i'}\int_{-\infty}^{\infty}\partial_x g_{i'}^{*}(x-y,t)
    \begin{pmatrix}
      \theta_1 \\
      \theta_2
    \end{pmatrix}
    (y,0)\, dy
  \end{aligned}
\end{align}
and
\begin{align}
  \label{eq:Ni}
  \begin{aligned}
    \mathcal{N}_i(x,t)
    & =\int_{0}^{t}\int_{-\infty}^{\infty}g_{ii}^{*}(x-y,t-s)\partial_x (n-n^*)(y,s)(y,s)\, dyds \\
    & \quad +\sum_{j=1}^{2}\int_{0}^{t}\int_{-\infty}^{\infty}(g_{ij}-g_{ij}^{*})(x-y,t-s)\partial_x n(y,s)\, dyds \\
    & \quad -\gamma_{i'}\int_{0}^{t}\int_{-\infty}^{\infty}\partial_x g_{i'i'}^{*}(x-y,t-s)\partial_x n^*(y,s)(y,s)\, dyds.
  \end{aligned}
\end{align}
Lemma~\ref{lem:vi_IE} may then be written as
\begin{equation}
  v_i(x,t)=\mathcal{I}_i(x,t)+\mathcal{N}_i(x,t).
\end{equation}
In the next two sections, we prove pointwise estimates for $\mathcal{I}_i(x,t)$ and $\mathcal{N}_i(x,t)$.

\subsubsection{Contribution from the initial data}
Our goal in this section is to prove the following pointwise estimates for $\mathcal{I}_i(x,t)$ defined by~\eqref{eq:Ii}.

\begin{lem}
  \label{lem:PWE_init}
  For any $n\geq 1$, there exist positive constants $\delta_n$ and $C_n$ such that if~\eqref{thm:main1:eq:smallness} holds, then we have
  \begin{equation}
    |\mathcal{I}_i(x,t)|\leq C\delta \Psi_{i;n}(x,t)
  \end{equation}
  for all $x\in \mathbb{R}$ and $t\geq 0$.
\end{lem}

\begin{proof}
  We assume $t\geq 1$ below (the case when $t<1$ is easier to handle). Let
  \begin{equation}
    \mathcal{I}_{i,1}(x,t)=\int_{-\infty}^{\infty}g_{i}^{*}(x-y,t)
    \begin{pmatrix}
      u_1-\theta_1 \\
      u_2-\theta_2
    \end{pmatrix}
    (y,0)\, dy
  \end{equation}
  and
  \begin{align}
    \mathcal{I}_{i,2}(x,t)
    & =\int_{-\infty}^{\infty}(g_i-g_{i}^{*})(x-y,t)
    \begin{pmatrix}
      u_1 \\
      u_2
    \end{pmatrix}
    (y,0)\, dy \\
    & \quad -\gamma_{i'}\int_{-\infty}^{\infty}\partial_x g_{i'}^{*}(x-y,t)
    \begin{pmatrix}
      \theta_1 \\
      \theta_2
    \end{pmatrix}
    (y,0)\, dy.
  \end{align}
  Then of course $\mathcal{I}_i(x,t)=\mathcal{I}_{i,1}(x,t)+\mathcal{I}_{i,2}(x,t)$.

  We first show
  \begin{equation}
    |\mathcal{I}_{i,1}(x,t)|\leq C\delta \Psi_{i;n}(x,t).
  \end{equation}
  For this purpose, set
  \begin{equation}
    \eta_j(x)=\int_{-\infty}^{x}(u_j-\theta_j)(y,0)\, dy
  \end{equation}
  and $\eta=(\eta_1 \, \eta_2)^{T}$. We then have
  \begin{equation}
    \mathcal{I}_{i,1}(x,t)=\int_{-\infty}^{\infty}g_{i}^{*}(x-y,t)\partial_x \eta(y)\, dy.
  \end{equation}
  By the definition of $M_i$, see~\eqref{eq:Mi}, we have
  \begin{equation}
    \eta_j(x)\coloneqq \int_{-\infty}^{x}(u_j-\theta_j)(y,0)\, dy=-\int_{x}^{\infty}(u_j-\theta_j)(y,0)\, dy.
  \end{equation}
  This and~\eqref{thm:main1:eq:smallness} imply
  \begin{equation}
    |\eta_j(x)|\leq C\delta (|x|+1)^{-\beta_n}.
  \end{equation}
  We first consider Case (i) $|x-\lambda_i t|\leq (t+1)^{1/2}$. In this case, integration by parts and~\eqref{eq:gij_aster} yield
  \begin{align}
    |\mathcal{I}_{i,1}(x,t)|
    & =\left| \int_{-\infty}^{\infty}\partial_x g_{i}^{*}(x-y,t)\eta(y)\, dy \right| \\
    & \leq C(t+1)^{-1}\int_{-\infty}^{\infty}|\eta_i(x)|\, dx\leq C\delta (t+1)^{-1}\leq C\delta \Psi_{i;n}(x,t).
  \end{align}
  We next consider Case (ii) $(t+1)^{1/2}<|x-\lambda_i t|<t+1$ with $x-\lambda_i t>0$ (the case when $x-\lambda_i t<0$ is similar). Again, by integration by parts,
  \begin{align}
    |\mathcal{I}_{i,1}(x,t)|
    & \leq C(t+1)^{-1}e^{-\frac{(x-\lambda_i t)^2}{Ct}}\int_{-\infty}^{(x-\lambda_i t)/2}|\eta_i(y)|\, dy \\
    & \quad +C\delta (t+1)^{-1}\int_{(x-\lambda_i t)/2}^{\infty}e^{-\frac{(x-y-\lambda_i t)^2}{Ct}}(y+1)^{-\beta_n}\, dy \\
    & \leq C\delta (t+1)^{-1}e^{-\frac{(x-\lambda_i t)^2}{Ct}}+C\delta (t+1)^{-1/2}(|x-\lambda_i t|+1)^{-\beta_n}\leq C\delta \Psi_{i;n}(x,t).
  \end{align}
  We finally consider Case (iii) $|x-\lambda_i t|\geq t+1$. For brevity, we assume $x-\lambda_i t>0$. In this case, by~\eqref{thm:main1:eq:smallness}, we have
  \begin{align}
    |\mathcal{I}_{i,1}(x,t)|
    & \leq C(t+1)^{-1/2}e^{-\frac{(x-\lambda_i t)^2}{Ct}}\int_{-\infty}^{(x-\lambda_i t)/2}|(u_i-\theta_i)(y,0)|\, dy \\
    & \quad +C\delta (t+1)^{-1/2}\int_{(x-\lambda_i t)/2}^{\infty}e^{-\frac{(x-y-\lambda_i t)^2}{2\nu t}}(y+1)^{-\alpha_n}\, dy \\
    & \leq C\delta e^{-\frac{t}{C}}e^{-\frac{(x-\lambda_i t)^2}{Ct}}+C\delta (|x-\lambda_i t|+1)^{-\alpha_n}\leq C\delta \Psi_{i;n}(x,t).
  \end{align}

  We next show
  \begin{equation}
    |\mathcal{I}_{i,2}(x,t)|\leq C\delta \Psi_{i;n}(x,t).
  \end{equation}
  Writing
  \begin{align}
    \mathcal{I}_{i,2}(x,t)
    & =\int_{-\infty}^{\infty}(g_i-g_{i}^{*}-\gamma_{i'}\partial_x g_{i'}^{*})(x-y,t)
    \begin{pmatrix}
      u_1 \\
      u_2
    \end{pmatrix}
    (y,0)\, dy \\
    & \quad +\gamma_{i'}\int_{-\infty}^{\infty}\partial_x g_{i'}^{*}(x-y,t)
    \begin{pmatrix}
      u_1-\theta_1 \\
      u_2-\theta_2
    \end{pmatrix}
    (y,0)\, dy
  \end{align}
  and applying~\eqref{eq:gi_PWE2}, we see that it suffices to show that
  \begin{align}
    A(x,t) & =\gamma_{i'}\partial_x \int_{-\infty}^{\infty}g_{i'}^{*}(x-y,t)
    \begin{pmatrix}
      u_1-\theta_1 \\
      u_2-\theta_2
    \end{pmatrix}
    (y,0)\, dy, \\
    B(x,t) & =(t+1)^{-1}\int_{-\infty}^{\infty}e^{-\frac{(x-y-\lambda_i t)^2}{Ct}}\left|
    \begin{pmatrix}
      u_1 \\
      u_2
    \end{pmatrix}
    \right| (y,0)\, dy, \\
    C(x,t) & =(t+1)^{-3/2}\int_{-\infty}^{\infty}e^{-\frac{(x-y-\lambda_{i'}t)^2}{Ct}}\left|
    \begin{pmatrix}
      u_1 \\
      u_2
    \end{pmatrix}
    \right| (y,0)\, dy \\
    D(x,t) & =e^{-\frac{c^2}{\nu}t}\left|
    \begin{pmatrix}
      u_1 \\
      u_2
    \end{pmatrix}
    \right| (x,0)
  \end{align}
  are all bounded by $C\delta \Psi_{i;n}(x,t)$. First, this is trivial for $D(x,t)$. Next, since $A(x,t)=\gamma_{i'}\partial_x \mathcal{I}_{i',1}(x,t)$, modifying the calculations above for $\mathcal{I}_{i,1}(x,t)$ yield the bound for $A(x,t)$. The bound for $B(x,t)$ is also obtained in a way similar to that for $\mathcal{I}_{i,1}(x,t)$ (except that we don't need $\eta_i$ in the analysis).

  Let us finally consider $C(x,t)$. First, Case (i) $|x-\lambda_{i'}t|\leq (t+1)^{1/2}$ is easy:
  \begin{equation}
    |C(x,t)|\leq C\delta (t+1)^{-3/2}\leq C\delta \Psi_{i;n}(x,t).
  \end{equation}
  Case (ii) $|x-\lambda_{i'}t|>(t+1)^{1/2}$ with $x-\lambda_{i'}t>0$ is as follows:
  \begin{align}
    |C(x,t)|
    & \leq C(t+1)^{-3/2}e^{-\frac{(x-\lambda_{i'}t)^2}{Ct}}\int_{-\infty}^{(x-\lambda_{i'}t)/2}\left|
    \begin{pmatrix}
      u_1 \\
      u_2
    \end{pmatrix}
    \right| (y,0)\, dy \\
    & \quad +C\delta (t+1)^{-3/2}\int_{(x-\lambda_{i'}t)/2}^{\infty}e^{-\frac{(x-y-\lambda_{i'}t)^2}{Ct}}(y+1)^{-\alpha_n}\, dy \\
    & \leq C\delta (t+1)^{-3/2}e^{-\frac{(x-\lambda_{i'}t)^2}{Ct}}+C\delta (t+1)^{-1}(|x-\lambda_{i'}t|+1)^{-\alpha_n}\leq C\delta \Psi_{i;n}(x,t).
  \end{align}
  The case when $x-\lambda_{i'}t<0$ is similar. This ends the proof of the lemma.
\end{proof}

\subsubsection{Contribution from the nonlinear terms}
Our goal in this section is to prove the following pointwise estimates for $\mathcal{N}_i(x,t)$ defined by~\eqref{eq:Ni}.

\begin{lem}
  \label{lem:PWE_nonlinear}
  For any $n\geq 1$, there exist positive constants $\delta_n$ and $C_n$ such that if~\eqref{thm:main1:eq:smallness} holds, then we have
  \begin{equation}
    |\mathcal{N}_i(x,t)|\leq C(\delta+P(t))^2 \Psi_{i;n}(x,t)
  \end{equation}
  for all $x\in \mathbb{R}$ and $t\geq 0$.
\end{lem}

To prove this lemma, we first prove some preparatory lemmas. To state these, we introduce the notation
\begin{equation}
  \mathcal{I}_i[f](x,t)\coloneqq \int_{0}^{t}\int_{-\infty}^{\infty}\partial_x \left\{ (t-s)^{-1/2}e^{-\frac{(x-y-\lambda_i(t-s))^2}{2\nu(t-s)}} \right\} f(y,s)\, dyds
\end{equation}
for a function $f=f(x,t)$.

\begin{lem}
  \label{lem:thetai'Xii}
  Let $n\geq 1$ and $\epsilon=\max(M_1,M_2)$. If $\epsilon$ is sufficiently small, we have
  \begin{equation}
    |\mathcal{I}_i[\theta_{i'}\Xi_i](x,t)|\leq C\epsilon^3 \Psi_{i;n}(x,t).
  \end{equation}
\end{lem}

\begin{proof}
  We only prove the lemma for $i=1$ (the other case is similar). By Lemma~\ref{lem:Xii}, we have
  \begin{equation}
    |(\theta_2 \Xi_1)(x,t)+(2c)^{-1}(\theta_{2}^{3}/2+\theta_{2}^{2}\Xi_2)(x,t)|\leq C\epsilon^3 \Theta_{4}(x,t;-c,\nu^*).
  \end{equation}
  Since
  \begin{equation}
    \Theta_4(x,t;-c,\nu^*)\leq C(t+1)^{-(2+1/2^{n+1})/2}\psi_n(x,t;-c),
  \end{equation}
  Lemma~\ref{lem:conv_psi_n_different_family} implies
  \begin{equation}
    |\mathcal{I}_1[\Theta_4(\cdot,\cdot;-c,\nu^*)](x,t)|\leq C\Psi_{1;n}(x,t).
  \end{equation}
  Next, note that~\eqref{eq:thetai} and~\eqref{eq:Xi} imply
  \begin{equation}
    L_2 \theta_2=-\partial_x(\theta_{2}^{2}/2), \quad L_2 \Xi_2=-\partial_x(\theta_{1}^{2}/2+\theta_1 \Xi_1+\theta_2 \Xi_2),
  \end{equation}
  where $L_2=\partial_t-c\partial_x-(\nu/2)\partial_{x}^{2}$. Using these, similar to the bound for $\zeta_{1;n}$ in the proof of Lemma~\ref{lem:xii_n}, we can show that
  \begin{equation}
    |\mathcal{I}_1[\theta_{2}^{3}/2+\theta_{2}^{2}\Xi_2](x,t)|\leq C\epsilon^3 \Psi_{1;n}(x,t).
  \end{equation}
  Combining these, we obtain the lemma.
\end{proof}

We next show the following.

\begin{lem}
  \label{lem:thetai'_square}
  Let $n\geq 1$ and $\epsilon=\max(M_1,M_2)$. If $\epsilon$ is sufficiently small, we have
  \begin{equation}
    |\mathcal{I}_i[\partial_x \theta_{i'}^{2}](x,t)|\leq C\epsilon^2 \Psi_{i;n}(x,t).
  \end{equation}
\end{lem}

\begin{proof}
  We only prove the lemma for $i=1$ (the other case is similar). Applying Lemma~\ref{lem:heat_equation_technique} yields
  \begin{equation}
    \mathcal{I}_1[\partial_x \theta_{2}^{2}](x,t)=(2c)^{-1}\sqrt{2\pi \nu}\partial_x \theta_{2}^{2}(x,t)+I_1(x,t)+I_2(x,t),
  \end{equation}
  where
  \begin{equation}
    I_1(x,t)=\int_{0}^{t^{1/2}}\int_{-\infty}^{\infty}\partial_x \left\{ (t-s)^{-1/2}e^{-\frac{(x-y-c(t-s))^2}{2\nu(t-s)}} \right\}\partial_x \theta_{2}^{2}(y,s)\, dyds
  \end{equation}
  and
  \begin{align}
    I_2(x,t)
    & =-(2c)^{-1}\int_{-\infty}^{\infty}(t-t^{1/2})^{-1/2}e^{-\frac{(x-y-c(t-\sqrt{t}))^2}{2\nu(t-\sqrt{t})}}\partial_x \theta_{2}^{2}(y,t^{1/2})\, dy \\
    & \quad -(2c)^{-1}\int_{t^{1/2}}^{t}\int_{-\infty}^{\infty}(t-s)^{-1/2}e^{-\frac{(x-y-c(t-s))^2}{2\nu(t-s)}}\partial_x L_2 \theta_{2}^{2}(y,s)\, dyds \\
    & \eqqcolon I_{21}(x,t)+I_{22}(x,t).
  \end{align}
  By Lemmas~\ref{lem:LZ97_I1} and~\ref{lem:LZ97_I21}, we obtain
  \begin{equation}
    |I_1(x,t)|+|I_{21}(x,t)|\leq C\epsilon^2 \log(t+2)\Theta_2(x,t;c,\nu^*)\leq C\epsilon^2 \Psi_{1;n}(x,t).
  \end{equation}
  For $I_{22}(x,t)$, we apply Lemma~\ref{lem:heat_equation_technique} (without the integral on $[0,t^{1/2}]$), which yields
  \begin{equation}
    \mathcal{I}_{22}(x,t)=-(2c)^{-2}\sqrt{2\pi \nu}L_2 \theta_{2}^{2}(x,t)+J_2(x,t),
  \end{equation}
  where
  \begin{align}
    J_2(x,t)
    & =(2c)^{-2}\int_{-\infty}^{\infty}(t-t^{1/2})^{-1/2}e^{-\frac{(x-y-c(t-\sqrt{t}))^2}{2\nu(t-\sqrt{t})}}L_2 \theta_{2}^{2}(y,t^{1/2})\, dy \\
    & \quad +(2c)^{-2}\int_{t^{1/2}}^{t}\int_{-\infty}^{\infty}(t-s)^{-1/2}e^{-\frac{(x-y-c(t-s))^2}{2\nu(t-s)}}L_{2}^{2}\theta_{2}^{2}(y,s)\, dyds \\
    & \eqqcolon J_{21}(x,t)+J_{22}(x,t).
  \end{align}
  By some tedious calculations, we obtain
  \begin{equation}
    L_2 \theta_{2}^{2}=-2\partial_x(\theta_{2}^{3}/3)-\nu(\partial_x \theta_2)^2
  \end{equation}
  and
  \begin{equation}
    L_{2}^{2}\theta_{2}^{2}=\partial_{x}^{2}(\theta_{2}^{4}/2)+\nu \partial_x[(\partial_x \theta_2)\partial_x \theta_{2}^{2}]+\nu (\partial_x \theta_2)\partial_{x}^{2}\theta_{2}^{2}+\nu^2(\partial_{x}^{2}\theta_2)^2.
  \end{equation}
  Since $|L_2 \theta_{2}^{2}(x,t)|\leq C\epsilon^2 \Theta_4(x,t;-c,\nu^*)$, Lemma~\ref{lem:LZ97_I21} yields
  \begin{equation}
    |J_{21}(x,t)|\leq C\epsilon^2 \Theta_{5/2}(x,t;c,\nu^*)\leq C\epsilon^2 \Psi_{1;n}(x,t).
  \end{equation}
  And since $|L_{2}^{2}\theta_{2}^{2}(x,t)|\leq C\epsilon^2 \Theta_6(x,t;-c,\nu^*)$, Lemma~\ref{lem:LZ97_I22} implies
  \begin{equation}
    |J_{22}(x,t)|\leq C\epsilon^2 [(x-c(t+1))^2+(t+1)]^{-5/4}\leq C\epsilon^2 \Psi_{1;n}(x,t;c).
  \end{equation}
  This proves the lemma.
\end{proof}

We can similarly show the following.

\begin{lem}
  \label{lem:thetai'Xii'}
  Let $n\geq 1$ and $\epsilon=\max(M_1,M_2)$. If $\epsilon$ is sufficiently small, we have
  \begin{equation}
    |\mathcal{I}_i[\partial_x(\theta_{i'}\Xi_{i'})](x,t)|\leq C\epsilon^3 \Psi_{i;n}(x,t).
  \end{equation}
\end{lem}

We move on to prove the following.

\begin{lem}
  \label{lem:thetai'vi'}
  Let $n\geq1$. If $\delta$ defined by~\eqref{thm:main1:eq:smallness} is sufficiently small, we have
  \begin{equation}
    |\mathcal{I}_i[\theta_{i'}v_{i'}](x,t)|\leq C(\delta+P(t))^2 \Psi_{i;n}(x,t).
  \end{equation}
  Here $P(t)$ is defined by~\eqref{eq:P}.
\end{lem}

\begin{proof}
  We only prove the lemma for $i=1$ (the other case is similar). Set $f=\theta_2 v_2$. Then Lemma~\ref{lem:heat_equation_technique} implies
  \begin{equation}
    \mathcal{I}_1[\theta_2 v_2](x,t)=(2c)^{-1}\sqrt{2\pi \nu}f(x,t)+I_1(x,t)+I_2(x,t),
  \end{equation}
  where
  \begin{equation}
    I_1(x,t)=\int_{0}^{t^{1/2}}\int_{-\infty}^{\infty}\partial_x \left\{ (t-s)^{-1/2}e^{-\frac{(x-y-c(t-s))^2}{2\nu(t-s)}} \right\} f(y,s)\, dyds
  \end{equation}
  and
  \begin{align}
    I_2(x,t)
    & =-(2c)^{-1}\int_{-\infty}^{\infty}(t-t^{1/2})^{-1/2}e^{-\frac{(x-y-c(t-\sqrt{t}))^2}{2\nu(t-\sqrt{t})}}f(y,t^{1/2})\, dy \\
    & \quad -(2c)^{-1}\int_{t^{1/2}}^{t}\int_{-\infty}^{\infty}(t-s)^{-1/2}e^{-\frac{(x-y-c(t-s))^2}{2\nu(t-s)}}L_2 f(y,s)\, dyds \\
    & \eqqcolon I_{21}(x,t)+I_{22}(x,t).
  \end{align}
  By Lemmas~\ref{lem:LZ97_I1} and~\ref{lem:LZ97_I21}, we obtain
  \begin{equation}
    |I_1(x,t)|+|I_{21}(x,t)|\leq C\delta P(t)\Theta_{\alpha_{n+1}}(x,t;c,\nu^*)\leq C\delta P(t)\Psi_{1;n}(x,t).
  \end{equation}
  To bound $I_{22}(x,t)$, we first note that
  \begin{align}
    L_2 v_2
    & =L_2(u_2-\theta_2-\Xi_2-\gamma_1 \partial_x \theta_1-\gamma_1 \partial_x \Xi_1) \\
    & =\frac{\nu}{2}\partial_{x}^{2}u_1+\partial_x n+\partial_x(\theta_{2}^{2}/2)+\partial_x(\theta_{1}^{2}/2+\theta_1 \Xi_1+\theta_2 \Xi_2) \\
    & \quad +\gamma_1 \partial_{x}^{2}(\theta_{1}^{2}/2+2c\theta_1)+\gamma_1 \partial_{x}^{2}(\theta_{2}^{2}/2+\theta_1 \Xi_1+\theta_2 \Xi_2+2c\Xi_1) \\
    & =\frac{\nu}{2}\partial_{x}^{2}(u_1-\theta_1-\Xi_1)+\partial_x(n-n_*)-\gamma_1 \partial_{x}^{2}n_*.
  \end{align}
  Then set
  \begin{equation}
    F=\frac{\nu}{2}\partial_x(u_1-\theta_1-\Xi_1)+n-n_*-\gamma_1 \partial_x n_*
  \end{equation}
  and
  \begin{equation}
    G=\theta_2 F-\nu v_2 \partial_x \theta_2.
  \end{equation}
  We note that
  \begin{equation}
    L_2 f=-v_2 \partial_x(\theta_{2}^{2}/2)+\partial_x G-(\partial_x \theta_2)F+\nu v_2 \partial_{x}^{2}\theta_2.
  \end{equation}
  By~\cite[Theorem~2.6 and Remark~2.8]{LZ97}, we have\footnote{The decay estimate for $\partial_{x}^{2}u$ is not explicitly stated in the theorem but is shown in its proof (see~\cite[p.~107]{LZ97}). Note that this is where the $H^6$-regularity is used.}
  \begin{equation}
    |\partial_x(u_1-\theta_1)(x,t)|\leq C\delta(t+1)^{-1/2}\Psi_{1;0}(x,t), \quad |\partial_{x}^{2}u_1(x,t)|\leq C\delta(t+1)^{-3/2}.
  \end{equation}
  In addition, applying Taylor's theorem, we see that
  \begin{equation}
    |(n-n_*)(x,t)|\leq C(\delta+P(t))^2(t+1)^{-1/2}[\psi_n(x,t;c)+\psi_n(x,t;-c)].
  \end{equation}
  These imply
  \begin{equation}
    |G(x,t)|\leq C(\delta+P(t))^2 (t+1)^{-1}\Theta_{\alpha_n}(x,t;-c,\nu^*)
  \end{equation}
  and
  \begin{equation}
    |L_2 f(x,t)-\partial_x G(x,t)|\leq C(\delta+P(t))^2 (t+1)^{-3/2}\Theta_{\alpha_n}(x,t;-c,\nu^*).
  \end{equation}
  Using these and integration by parts, we get
  \begin{align}
    |I_{22}(x,t)|
    & \leq C(\delta+P(t))^2 \int_{t^{1/2}}^{t/2}\int_{-\infty}^{\infty}(t-s)^{-1/2}e^{-\frac{(x-y-c(t-s))^2}{C(t-s)}}(s+1)^{-3/2}\psi_n(y,s;-c)\, dyds \\
    & \quad +C(\delta+P(t))^2 \int_{t^{1/2}}^{t}\int_{-\infty}^{\infty}(t-s)^{-1}e^{-\frac{(x-y-c(t-s))^2}{C(t-s)}}(s+1)^{-1}\psi_n(y,s;-c)\, dyds.
  \end{align}
  Applying Lemmas~\ref{lem:conv_psi_n_different_family} and~\ref{lem:conv_psi_n_different_family_sqrt_t_half}, we obtain
  \begin{equation}
    |I_{22}(x,t)|\leq C(\delta+P(t))^2 \Psi_{1;n}(x,t).
  \end{equation}
  This proves the lemma.
\end{proof}

The lemma below can be shown in a similar manner.

\begin{lem}
  \label{lem:Xii'_square}
  Let $n\geq 1$. If $\delta$ defined by~\eqref{thm:main1:eq:smallness} is sufficiently small, we have
  \begin{equation}
    |\mathcal{I}_i[\Xi_{i'}^{2}](x,t)|+|\mathcal{I}_i[\Xi_{i'}v_{i'}](x,t)|\leq C(\delta+P(t))^3 \Psi_{i;n}(x,t).
  \end{equation}
\end{lem}

Set
\begin{equation}
  n_a=-\frac{p''(1)^2}{8c^2}(v-1)^2, \quad n_b=-\frac{\nu p''(1)}{4c^2}(v-1)u_x, \quad n_c=n-n_a-n_b.
\end{equation}
Of course $n=n_a+n_n+n_c$. Correspondingly, set
\begin{align}
  \begin{aligned}
    \mathcal{N}_{i,a}(x,t)
    & =\int_{0}^{t}\int_{-\infty}^{\infty}g_{ii}^{*}(x-y,t-s)\partial_x (n_a-n^*)(y,s)(y,s)\, dyds \\
    & \quad +\sum_{j=1}^{2}\int_{0}^{t}\int_{-\infty}^{\infty}(g_{ij}-g_{ij}^{*})(x-y,t-s)\partial_x n_a(y,s)\, dyds \\
    & \quad -\gamma_{i'}\int_{0}^{t}\int_{-\infty}^{\infty}\partial_x g_{i'i'}^{*}(x-y,t-s)\partial_x n^*(y,s)(y,s)\, dyds,
  \end{aligned}
\end{align}
\begin{align}
  \begin{aligned}
    \mathcal{N}_{i,b}(x,t)
    & =\int_{0}^{t}\int_{-\infty}^{\infty}g_{ii}^{*}(x-y,t-s)\partial_x n_b(y,s)(y,s)\, dyds \\
    & \quad +\sum_{j=1}^{2}\int_{0}^{t}\int_{-\infty}^{\infty}(g_{ij}-g_{ij}^{*})(x-y,t-s)\partial_x n_b(y,s)\, dyds,
  \end{aligned}
\end{align}
and
\begin{align}
  \begin{aligned}
    \mathcal{N}_{i,c}(x,t)
    & =\int_{0}^{t}\int_{-\infty}^{\infty}g_{ii}^{*}(x-y,t-s)\partial_x n_c(y,s)(y,s)\, dyds \\
    & \quad +\sum_{j=1}^{2}\int_{0}^{t}\int_{-\infty}^{\infty}(g_{ij}-g_{ij}^{*})(x-y,t-s)\partial_x n_c(y,s)\, dyds.
  \end{aligned}
\end{align}
Then $\mathcal{N}_i(x,t)=\mathcal{N}_{i,a}(x,t)+\mathcal{N}_{i,b}(x,t)+\mathcal{N}_{i,c}(x,t)$; see~\eqref{eq:Ni}.

We next prove the following.

\begin{lem}
  Let $n\geq 1$. If $\delta$ defined by~\eqref{thm:main1:eq:smallness} is sufficiently small, we have
  \begin{equation}
    \label{lem:PWE_nonlinear:Na}
    |\mathcal{N}_{i,a}(x,t)|\leq C(\delta+P(t))^2 \Psi_{i;n}(x,t).
  \end{equation}
\end{lem}

\begin{proof}
  Let $i=1$ (the case of $i=2$ is similar). By integration by parts, we have
  \begin{align}
    \mathcal{N}_{1,a}(x,t)
    & =\int_{0}^{t}\int_{-\infty}^{\infty}\partial_x g_{11}^{*}(x-y,t-s)(n_a-n^*)(y,s)\, dyds \\
    & \quad +\sum_{j=1}^{2}\int_{0}^{t}\int_{-\infty}^{\infty}\partial_x (g_{1j}-g_{1j}^{*})(x-y,t-s)n_a(y,s)\, dyds \\
    & \quad -\gamma_2 \int_{0}^{t}\int_{-\infty}^{\infty}\partial_{x}^{2}g_{22}^{*}(x-y,t-s)n^*(y,s)\, dyds.
  \end{align}
  By some tedious calculations, we can show that
  \begin{align}
    \begin{aligned}
      & |(n_a-n^*)(x,t)-[\theta_2 \Xi_1+\gamma_2 \partial_x(\theta_{2}^{2}/2)+\gamma_2 \partial_x(\theta_2 \Xi_2)-\theta_2 v_2-\Xi_{2}^{2}/2-\Xi_2 v_2](x,t)| \\
      & \leq C(\delta+P(t))^2[(t+1)^{-1/2}\psi_n(x,t;c)+(t+1)^{-\alpha_n/2}\psi_n(x,t;-c)].
    \end{aligned}
  \end{align}
  Then Lemmas~\ref{lem:thetai'Xii}--\ref{lem:Xii'_square},~\ref{lem:conv_psi_n_same_family}, and~\ref{lem:conv_psi_n_different_family} yield
  \begin{align}
    \left| \int_{0}^{t}\int_{-\infty}^{\infty}\partial_x g_{11}^{*}(x-y,t-s)(n_a-n^*)(y,s)\, dyds \right| \leq C(\delta+P(t))^2 \Psi_{i;n}(x,t).
  \end{align}
  It remains to show that
  \begin{equation}
    |I(x,t)|\leq C(\delta+P(t))^2 \Psi_{1;n}(x,t),
  \end{equation}
  where
  \begin{align}
    I(x,t)
    & =\sum_{j=1}^{2}\int_{0}^{t}\int_{-\infty}^{\infty}\partial_x (g_{1j}-g_{1j}^{*})(x-y,t-s)n_a(y,s)\, dyds \\
    & \quad -\gamma_2 \int_{0}^{t}\int_{-\infty}^{\infty}\partial_{x}^{2}g_{22}^{*}(x-y,t-s)n^*(y,s)\, dyds.
  \end{align}
  We define the decomposition $I(x,t)=I_1(x,t)+I_2(x,t)$ by
  \begin{equation}
    I_1(x,t)=\sum_{j=1}^{2}\int_{0}^{t}\int_{-\infty}^{\infty}\partial_x (g_{1j}-g_{1j}^{*})(x-y,t-s)(n_a-n^*)(y,s)\, dyds
  \end{equation}
  and
  \begin{align}
    I_2(x,t)=
    & \sum_{j=1}^{2}\int_{0}^{t}\int_{-\infty}^{\infty}\partial_x (g_{1j}-g_{1j}^{*})(x-y,t-s)n^*(y,s)\, dyds \\
    & \quad -\gamma_2 \int_{0}^{t}\int_{-\infty}^{\infty}\partial_{x}^{2}g_{22}^{*}(x-y,t-s)n^*(y,s)\, dyds.
  \end{align}

  We first consider $I_1(x,t)$. By~\eqref{eq:gi_PWE1}, to show that $|I_1(x,t)|$ is bounded by $C(\delta+P(t))^2 \Psi_{1;n}(x,t)$, it suffices to prove the same bound for
  \begin{equation}
    \int_{0}^{t}\int_{-\infty}^{\infty}(t-s)^{-1}(t-s+1)^{-1/2}e^{-\frac{(x-y-\lambda_j(t-s))^2}{C(t-s)}}|(n_a-n^*)(y,s)|\, dyds \quad (j=1,2)
  \end{equation}
  and
  \begin{equation}
    \int_{0}^{t}e^{-\frac{c^2}{\nu}(t-s)}|(n_a-n^*)(x,s)|\, ds.
  \end{equation}
  The term corresponding to $\delta^{(1)}(x)$ is not needed since $q_{10}=(1/2 \, -1/2)^{T}$. Noting that
  \begin{equation}
    |(n_a-n^*)(x,t)|\leq C(\delta+P(t))^2(t+1)^{-1/2}[\psi_n(x,t;c)+\psi_n(x,t;-c)],
  \end{equation}
  Lemmas~\ref{lem:conv_psi_n_same_family},~\ref{lem:conv_psi_n_different_family}, and~\ref{lem:conv_exp} imply the desired bounds for the two integrals above.

  We next consider $I_2(x,t)$. We have $I_2(x,t)=I_{21}(x,t)+I_{22}(x,t)$ with
  \begin{align}
    I_{21}(x,t)
    & =-\sum_{j=1}^{2}\int_{0}^{t}\int_{-\infty}^{\infty}\partial_x (g_{1j}-g_{1j}^{*})(x-y,t-s)(\theta_{1}^{2}/2+\theta_1 \Xi_1)(y,s)\, dyds \\
    & \quad +\gamma_2 \int_{0}^{t}\int_{-\infty}^{\infty}\partial_{x}^{2}g_{22}^{*}(x-y,t-s)(\theta_{1}^{2}/2+\theta_1 \Xi_1)(y,s)\, dyds
  \end{align}
  and
  \begin{align}
    I_{22}(x,t)
    & =-\sum_{j=1}^{2}\int_{0}^{t}\int_{-\infty}^{\infty}\partial_x (g_{1j}-g_{1j}^{*})(x-y,t-s)(\theta_{2}^{2}/2+\theta_2 \Xi_2)\, dyds \\
    & \quad +\gamma_2 \int_{0}^{t}\int_{-\infty}^{\infty}\partial_{x}^{2}g_{22}^{*}(x-y,t-s)(\theta_{2}^{2}/2+\theta_2 \Xi_2)\, dyds.
  \end{align}
  Taking into account~\eqref{eq:gi_PWE1} and~\eqref{eq:gi_PWE2}, Lemmas~\ref{lem:conv_psi_n_same_family},~\ref{lem:conv_psi_n_different_family}, and~\ref{lem:conv_exp} yield $|I_{21}(x,t)|\leq C(\delta+P(t))^2 \Psi_{1;n}(x,t)$ (we divide the domain of temporal integration into $[0,t/2]$ and $[t/2,t]$ then use integration by parts before applying the lemmas). For $I_{22}(x,t)$, we proceed as follows. Using the technique in the proof of Lemma~\ref{lem:heat_equation_technique}, we obtain
  \begin{align}
    I_{22}(x,t)
    & =-\sum_{j=1}^{2}\int_{0}^{t^{1/2}}\int_{-\infty}^{\infty}\partial_x(g_{1j}-g_{1j}^{*})(x-y,t-s)(\theta_{2}^{2}/2+\theta_2 \Xi_2)(y,s)\, dyds \\
    & \quad +\gamma_2 \int_{0}^{t^{1/2}}\int_{-\infty}^{\infty}\partial_{x}^{2}g_{22}^{*}(x-y,t-s)(\theta_{2}^{2}/2+\theta_2 \Xi_2)(y,s)\, dyds \\
    & \quad -\frac{1}{2c}\sum_{j=1}^{2}\int_{t^{1/2}}^{t}\int_{-\infty}^{\infty}L_1(g_{1j}-g_{1j}^{*})(x-y,t-s)(\theta_{2}^{2}/2+\theta_2 \Xi_2)(y,s)\, dyds \\
    & \quad +\frac{\nu}{4c}\int_{t^{1/2}}^{t}\int_{-\infty}^{\infty}\partial_{x}^{2}g_{22}^{*}(x-y,t-s)(\theta_{2}^{2}/2+\theta_2 \Xi_2)(y,s)\, dyds \\
    & \quad +\frac{1}{2c}\sum_{j=1}^{2}\int_{t^{1/2}}^{t}\int_{-\infty}^{\infty}(g_{1j}-g_{1j}^{*})(x-y,t-s)L_2(\theta_{2}^{2}/2+\theta_2 \Xi_2)(y,s)\, dyds \\
    & \quad +\frac{1}{2c}\sum_{j=1}^{2}\int_{-\infty}^{\infty}(g_{1j}-g_{1j}^{*})(x-y,t-t^{1/2})(\theta_{2}^{2}/2+\theta_2 \Xi_2)(y,t^{1/2})\, dyds,
  \end{align}
  where $L_i=\partial_t+\lambda_i \partial_x-(\nu/2)\partial_{x}^{2}$. Here we used $\lim_{t\to 0}(g_{1j}-g_{1j}^{*})(x,t)=0$. The sum of the first two terms on the right-hand side can be bounded using Lemmas~\ref{lem:conv_psi_n_same_family},~\ref{lem:conv_psi_n_different_family_0_sqrt}, and~\ref{lem:conv_exp}; the sum of the third and the fourth term can be bounded using Lemmas~\ref{lem:conv_psi_n_same_family},~\ref{lem:conv_psi_n_different_family},~\ref{lem:conv_exp}, and the relation
  \begin{equation}
    L_1(g_{1j}-g_{1j}^{*})=(\nu/2)\partial_{x}^{2}g_{2j}.
  \end{equation}
  To bound the sum of the fifth and the sixth term, noting that
  \begin{equation}
    |L_2(\theta_{2}^{2}/2+\theta_2 \Xi_2)(x,t)|\leq C\delta^2 \Theta_4(x,t;-c,\nu^*),
  \end{equation}
  it suffices to show that the following integrals are bounded by $C(\delta+P(t))^2 \Psi_{1;n}(x,t)$:
  \begin{equation}
    A(x,t)=\int_{t^{1/2}}^{t}\int_{-\infty}^{\infty}(t-s)^{-1/2}(t-s+1)^{-1/2}e^{-\frac{(x-y-c(t-s))^2}{C(t-s)}}\Theta_4(y,s;-c,\nu^*)\, dyds,
  \end{equation}
  \begin{equation}
    B(x,t)=\int_{-\infty}^{\infty}(t-t^{1/2})^{-1}e^{-\frac{(x-y-c(t-\sqrt{t}))^2}{C(t-\sqrt{t})}}\Theta_2(y,t^{1/2};-c,\nu^*)\, dy,
  \end{equation}
  \begin{equation}
    C(x,t)=\int_{t^{1/2}}^{t}\int_{-\infty}^{\infty}(t-s)^{-1}(t-s+1)^{-1/2}e^{-\frac{(x-y+c(t-s))^2}{C(t-s)}}\Theta_4(y,s;-c,\nu^*)\, dyds,
  \end{equation}
  \begin{equation}
    D(x,t)=\int_{-\infty}^{\infty}(t-t^{1/2})^{-1}(t-t^{1/2}+1)^{-1/2}e^{-\frac{(x-y+c(t-\sqrt{t}))^2}{C(t-\sqrt{t})}}\Theta_2(y,t^{1/2};-c,\nu^*)\, dy,
  \end{equation}
  and
  \begin{align}
    E(x,t)
    & =\int_{t^{1/2}}^{t}\int_{-\infty}^{\infty}\partial_x g_{22}^{*}(x-y,t-s)L_2(\theta_{2}^{2}/2+\theta_2 \Xi_2)(y,s)\, dyds \\
    & \quad +\int_{-\infty}^{\infty}\partial_x g_{22}^{*}(x-y,t-t^{1/2})(\theta_{2}^{2}/2+\theta_2 \Xi_2)(y,t^{1/2})\, dy.
  \end{align}
  We can bound $A(x,t)$ using Lemma~\ref{lem:LZ97_I22}, $B(x,t)$ and $D(x,t)$ by Lemma~\ref{lem:LZ97_I21}, and $C(x,t)$ by Lemma~\ref{lem:conv_psi_n_same_family}. Finally, we consider $E(x,t)$. Taking into account $L_2 g_{22}^{*}=0$ and $\lim_{t\to 0}\partial_x g_{22}^{*}(x-y,t)=\delta^{(1)}(x-y)$, integration by parts applied to the operator $L_2$ yields
  \begin{equation}
    E(x,t)=\int_{-\infty}^{\infty}\delta^{(1)}(x-y)(\theta_{2}^{2}/2+\theta_2 \Xi_2)(y,t)\, dy=\partial_x(\theta_{2}^{2}/2+\theta_2 \Xi_2)(x,t).
  \end{equation}
  Hence $|E(x,t)|\leq C(\delta+P(t))^2 \Psi_{1;n}(x,t)$. This ends the proof of the lemma.
\end{proof}

The lemma below can be proved similarly (and in fact much simply).

\begin{lem}
  Let $n\geq 1$. If $\delta$ defined by~\eqref{thm:main1:eq:smallness} is sufficiently small, we have
  \begin{equation}
    \label{lem:PWE_nonlinear:Nb}
    |\mathcal{N}_{i,b}(x,t)|+|\mathcal{N}_{i,c}(x,t)|\leq C(\delta+P(t))^2 \Psi_{i;n}(x,t).
  \end{equation}
\end{lem}

Combining Lemmas~\ref{lem:PWE_nonlinear:Na} and~\ref{lem:PWE_nonlinear:Nb}, the proof of Lemma~\ref{lem:PWE_nonlinear} is complete.

\subsubsection{Final step of the proof}
\label{sec:final}
The remaining step of the proof is standard. By Lemma~\ref{lem:PWE_init} and~\ref{lem:PWE_nonlinear}, we obtain
\begin{equation}
  \label{eq:P_ineq:var}
  P(t)\leq C\delta+C(\delta+P(t))^2 \leq C_1 \delta+C_2 P(t)^2
\end{equation}
for some $C_1,C_2>0$. Here $P(t)$ is defined by~\eqref{eq:P}. When $\delta$ is sufficiently small, the line $y=p$ and the parabola $y=C_1 \delta+C_2 p^2$ intersect at $p=p_1$ and $p_2$, where
\begin{equation}
  p_1=\frac{1-\sqrt{1-4C_1 C_2 \delta}}{2C_2}, \quad p_2=\frac{1+\sqrt{1-4C_1 C_2 \delta}}{2C_2}.
\end{equation}
Note that $C_1 \delta \leq p_1<p_2$. By~\eqref{eq:P_ineq:var}, we either have $P(t)\leq p_1$ or $P(t)\geq p_2$. Since $P(t)$ is continuous in $t$, if $P(0)\leq p_1$, then $P(t)\leq p_1$ for all $t\geq 0$. By~\eqref{thm:main1:eq:smallness}, taking $C_1$ sufficiently large, we indeed have $P(0)\leq C_1 \delta \leq p_1$. Therefore, we conclude that
\begin{equation}
  P(t)\leq p_1 \leq \frac{1-(1-4C_1 C_2 \delta)}{2C_2}\leq 2C_1 \delta.
\end{equation}
This ends the proof of Theorem~\ref{thm:main1}.

\appendix
\renewcommand*{\thesection}{\Alph{section}}
\section{Lemmas on convolutions involving a heat kernel}
\begin{lem}
  \label{lem:heat_equation_technique}
  Suppose that $f=f(x,t)$ is a $C^2$-smooth function on $\mathbb{R}\times (0,\infty)$. Let $\lambda \neq \lambda'$ and $\nu>0$, and set $L_{\lambda'}=\partial_t+\lambda' \partial_x-(\nu/2)\partial_{x}^{2}$. Then for $t\geq 1$, the function $I(x,t)$ defined by
  \begin{equation}
    I(x,t)=\int_{0}^{t}\int_{-\infty}^{\infty}(t-s)^{-1/2}e^{-\frac{(x-y-\lambda(t-s))^2}{2\nu(t-s)}}\partial_x f(y,s)\, dyds
  \end{equation}
  can be written as
  \begin{equation}
    I(x,t)=(\lambda-\lambda')^{-1}\sqrt{2\pi \nu}f(x,t)+I_1(x,t)+I_2(x,t),
  \end{equation}
  where
  \begin{equation}
    I_1(x,t)=\int_{0}^{t^{1/2}}\int_{-\infty}^{\infty}(t-s)^{-1/2}e^{-\frac{(x-y-\lambda(t-s))^2}{2\nu(t-s)}}\partial_x f(y,s)\, dyds,
  \end{equation}
  and
  \begin{align}
    I_2(x,t)
    & =-(\lambda-\lambda')^{-1}\int_{-\infty}^{\infty}(t-t^{1/2})^{-1/2}e^{-\frac{(x-y-\lambda(t-\sqrt{t}))^2}{2\nu(t-\sqrt{t})}}f(y,t^{1/2})\, dy \\
    & \quad -(\lambda-\lambda')^{-1}\int_{t^{1/2}}^{t}\int_{-\infty}^{\infty}(t-s)^{-1/2}e^{-\frac{(x-y-\lambda(t-s))^2}{2\nu(t-s)}}L_{\lambda'}f(y,s)\, dyds.
  \end{align}
\end{lem}

\begin{proof}
  Let
  \begin{equation}
    g_{\lambda}(x,t)=t^{-1/2}e^{-\frac{(x-\lambda t)^2}{2\nu t}}.
  \end{equation}
  Dividing the domain of temporal integration, we get
  \begin{align}
    I(x,t)
    & =\int_{0}^{t^{1/2}}\int_{-\infty}^{\infty}g_{\lambda}(x-y,t-s)\partial_x f(y,s)\, dyds \\
    & \quad +\int_{t^{1/2}}^{t}\int_{-\infty}^{\infty}g_{\lambda}(x-y,t-s)\partial_x f(y,s)\, dyds.
  \end{align}
  The first term on the right-hand side is $I_1(x,t)$. For the second term, integration by parts yields
  \begin{align}
    & \int_{t^{1/2}}^{t}\int_{-\infty}^{\infty}g_{\lambda}(x-y,t-s)\partial_x f(y,s)\, dyds \\
    & =-(\lambda-\lambda')^{-1}\int_{t^{1/2}}^{t}\int_{-\infty}^{\infty}g_{\lambda}(x-y,t-s) \\
    & \qquad \qquad \cdot [-\partial_s+\partial_s-\lambda \partial_y+\lambda' \partial_y+(\nu/2)\partial_{y}^{2}-(\nu/2)\partial_{y}^{2}]f(y,s)\, dyds \\
    & =(\lambda-\lambda')^{-1}\int_{t^{1/2}}^{t}\int_{-\infty}^{\infty}L_{\lambda}g_{\lambda}(x-y,t-s)f(y,s)\, dyds \\
    & \quad -(\lambda-\lambda')^{-1}\int_{t^{1/2}}^{t}\int_{-\infty}^{\infty}g_{\lambda}(x-y,t-s)L_{\lambda'}f(y,s)\, dyds \\
    & \quad +(\lambda-\lambda')^{-1}\sqrt{2\pi \nu}f(x,t)-(\lambda-\lambda')^{-1}\int_{-\infty}^{\infty}g_{\lambda}(x-y,t-t^{1/2})f(x,t^{1/2})\, dyds,
  \end{align}
  where $L_{\lambda}=\partial_s+\lambda \partial_y-(\nu/2)\partial_{y}^{2}$. Here we used $\lim_{s\to t}g_{\lambda}(x-y,t-s)=\sqrt{2\pi \nu}\delta(x-y)$. The lemma follows from the equality above by noting that $L_{\lambda}g_{\lambda}=0$.
\end{proof}

In the lemmas below, $C$ and $\nu^*$ denote generic large constants. We remind the reader that $\Theta_{\alpha}(x,t;\lambda,\mu)$ is defined by~\eqref{eq:Theta}.

\begin{lem}
  \label{lem:LZ97_I1}
  Let $\lambda \neq \lambda'$, $\mu>0$, and $0<\alpha \leq 3$. Then we have
  \begin{align}
    \label{lem:LZ97_I1_ineq1}
    \begin{aligned}
      & \int_{0}^{t^{1/2}}\int_{-\infty}^{\infty}(t-s)^{-1}e^{-\frac{(x-y-\lambda(t-s))^2}{\mu(t-s)}}\Theta_{\alpha}(y,s;\lambda',\mu)\, dyds \\
      & \leq
      \begin{dcases}
        C\Theta_{(\alpha+1)/2}(x,t;\lambda,\nu^*) & \text{if $\alpha \neq 3$}, \\
        C\log(t+2)\Theta_{(\alpha+1)/2}(x,t;\lambda,\nu^*) & \text{if $\alpha=3$}.
      \end{dcases}
    \end{aligned}
  \end{align}
\end{lem}

\begin{proof}
  See the analysis of $I_1(x,t)$ in the proof of~\cite[Lemma~3.4]{LZ97}.
\end{proof}

\begin{lem}
  \label{lem:LZ97_I21}
  Let $\lambda,\lambda' \in \mathbb{R}$, $\mu>0$, and $\alpha>0$ (not necesarily $\lambda \neq \lambda'$). Then we have
  \begin{equation}
    \label{lem:LZ97_I21_ineq1}
    \int_{-\infty}^{\infty}e^{-\frac{(x-y-\lambda(t-\sqrt{t}))^2}{\mu(t-\sqrt{t})}}\Theta_{\alpha}(y,t^{1/2};\lambda',\mu)\, dyds\leq C\Theta_{(\alpha-1)/2}(x,t;\lambda,\nu^*).
  \end{equation}
\end{lem}

\begin{proof}
  See the analysis of $I_{21}(x,t)$ in the proof of~\cite[Lemma~3.4]{LZ97}.
\end{proof}

\begin{lem}
  \label{lem:LZ97_I22}
  Let $\lambda \neq \lambda'$, $\mu>0$, and $\alpha>1$. Then we have
  \begin{align}
    \label{lem:LZ97_I22_ineq1}
    \begin{aligned}
      & \int_{t^{1/2}}^{t}\int_{-\infty}^{\infty}(t-s)^{-1/2}e^{-\frac{(x-y-\lambda(t-s))^2}{\mu(t-s)}}\Theta_{\alpha}(y,s;\lambda',\mu)\, dyds \\
      & \quad \leq C[(x-\lambda(t+1))^2+(t+1)]^{-(\alpha-1)/4}.
    \end{aligned}
  \end{align}
\end{lem}

\begin{proof}
  See the analysis of $I_{22}^{(1)}(x,t)$ in the proof of~\cite[Lemma~3.4]{LZ97}.
\end{proof}

\begin{lem}[{\cite[Lemma~3.2]{LZ97}}]
  Let $\lambda \in \mathbb{R}$, $\mu>0$, $\alpha \geq 0$, and $\beta>0$. Then we have
  \begin{align}
    \begin{aligned}
      & \int_{0}^{t/2}\int_{-\infty}^{\infty}(t-s)^{-1}(t+1-s)^{-\alpha/2}e^{-\frac{(x-y-\lambda(t-s))^2}{\mu(t-s)}}\Theta_{\beta}(y,s;\lambda,\mu)\, dyds \\
      & \leq
      \begin{dcases}
        C\Theta_{\gamma}(x,t;\lambda,\nu^*) & \text{if $\beta \neq 3$}, \\
        C\log(t+2)\Theta_{\gamma}(x,t;\lambda,\nu^*) & \text{if $\beta=3$},
      \end{dcases}
    \end{aligned}
  \end{align}
  where $\gamma=\alpha+\min(\beta,3)-1$. We also have
  \begin{align}
    \begin{aligned}
      & \int_{t/2}^{t}\int_{-\infty}^{\infty}(t-s)^{-1}(t+1-s)^{-\alpha/2}e^{-\frac{(x-y-\lambda(t-s))^2}{\mu(t-s)}}\Theta_{\beta}(y,s;\lambda,\mu)\, dyds \\
      & \leq
      \begin{dcases}
        C\Theta_{\gamma}(x,t;\lambda,\nu^*) & \text{if $\alpha \neq 1$}, \\
        C\log(t+2)\Theta_{\gamma}(x,t;\lambda,\nu^*) & \text{if $\alpha=1$},
      \end{dcases}
    \end{aligned}
  \end{align}
  where $\gamma=\min(\alpha,1)+\beta-1$.
\end{lem}

We remind the reader that $\psi_n(x,t;\lambda)$ is defined in~\eqref{eq:aux}.

\begin{lem}
  \label{lem:conv_psi_n_same_family}
  Let $\lambda \in \mathbb{R}$, $\mu>0$, and $\alpha,\beta \geq 0$. Then we have
  \begin{align}
    \label{lem:conv_psi_n_same_family_ineq1}
    \begin{aligned}
      & \int_{0}^{t/2}\int_{-\infty}^{\infty}(t-s)^{-1}(t+1-s)^{-\alpha/2}e^{-\frac{(x-y-\lambda(t-s))^2}{\mu(t-s)}}(s+1)^{-\beta/2}\psi_n(y,s;\lambda)\, dyds \\
      & \leq
      \begin{dcases}
        C(t+1)^{-\gamma_1/2}\psi_n(x,t;\lambda) & \text{if $\beta \neq 3-\alpha_n$}, \\
        C\log(t+2)(t+1)^{-\gamma_1/2}\psi_n(x,t;\lambda) & \text{if $\beta=3-\alpha_n$},
      \end{dcases}
    \end{aligned}
  \end{align}
  where $\gamma_1=\alpha+\min(\beta,3-\alpha_n)-1$. We also have
  \begin{align}
    \label{lem:conv_psi_n_same_family_ineq2}
    \begin{aligned}
      & \int_{t/2}^{t}\int_{-\infty}^{\infty}(t-s)^{-1}(t+1-s)^{-\alpha/2}e^{-\frac{(x-y-\lambda(t-s))^2}{\mu(t-s)}}(s+1)^{-\beta/2}\psi_n(y,s;\lambda)\, dyds \\
      & \leq
      \begin{dcases}
        C(t+1)^{-\gamma_2/2}\psi_n(x,t;\lambda) & \text{if $\alpha \neq 1$}, \\
        C\log(t+2)(t+1)^{-\gamma_2/2}\psi_n(x,t;\lambda) & \text{if $\alpha=1$},
      \end{dcases}
    \end{aligned}
  \end{align}
  where $\gamma_2=\min(\alpha,1)+\beta-1$.
\end{lem}

\begin{proof}
  A straightforward (but lengthy) adaptation of the proof of~\cite[Lemma~A.7]{Koike22b} proves the lemma.
\end{proof}

For $\lambda,\lambda' \in \mathbb{R}$ and $K>0$, let
\begin{align}
  & \chi_K(x,t;\lambda,\lambda') \\
  & \coloneqq \mathrm{char}\left\{ \min(\lambda,\lambda')(t+1)+K\sqrt{t+1}\leq x\leq \max(\lambda,\lambda')(t+1)-K\sqrt{t+1} \right\},
\end{align}
where $\mathrm{char}\{ S \}$ is the indicator function of a set $S$.

\begin{lem}
  \label{lem:conv_psi_n_different_family}
  Let $\lambda \neq \lambda'$, $\mu>0$, $\alpha \geq 0$, and $0\leq \beta \leq 2\alpha_n$ ($\beta \neq 2$).\footnote{The case of $\beta=2$ is excluded just for simplicity.}~Then for $K>0$ large enough, we have
  \begin{align}
    \label{lem:conv_psi_n_different_family_ineq1}
    \begin{aligned}
      & \int_{0}^{t/2}\int_{-\infty}^{\infty}(t-s)^{-1}(t+1-s)^{-\alpha/2}e^{-\frac{(x-y-\lambda(t-s))^2}{\mu(t-s)}}(s+1)^{-\beta/2}\psi_n(y,s;\lambda')\, dyds \\
      & \leq
      \begin{dcases}
        C[(t+1)^{-\gamma_1/2}\psi_n(x,t;\lambda)+(t+1)^{-\gamma_{1}'/2}\psi_n(x,t;\lambda')] & \text{if $\beta \neq 3-\alpha_n$} \\
        C[\log(t+2)(t+1)^{-\gamma_1/2}\psi_n(x,t;\lambda)+(t+1)^{-\gamma_{1}'/2}\psi_n(x,t;\lambda')] & \text{if $\beta=3-\alpha_n$}
      \end{dcases} \\
      & \quad +C|x-\lambda(t+1)|^{-(\alpha_n+\min(\beta,\alpha_n+1)-1)/2}|x-\lambda'(t+1)|^{-\alpha/2-1/2}\chi_K(x,t;\lambda,\lambda'),
    \end{aligned}
  \end{align}
  where $\gamma_1=\alpha+\min(\beta,3-\alpha_n)-1$ and $\gamma_{1}'=\alpha+\min(\beta,2)-1$. We also have
  \begin{align}
    \label{lem:conv_psi_n_different_family_ineq2}
    \begin{aligned}
      & \int_{t/2}^{t}\int_{-\infty}^{\infty}(t-s)^{-1}(t+1-s)^{-\alpha/2}e^{-\frac{(x-y-\lambda(t-s))^2}{\mu(t-s)}}(s+1)^{-\beta/2}\psi_n(y,s;\lambda')\, dyds \\
      & \leq
      \begin{dcases}
        C(t+1)^{-\gamma_2/2}[\psi_n(x,t;\lambda)+\psi_n(x,t;\lambda')] & \text{if $\alpha \neq 1$} \\
        C\log(t+2)(t+1)^{-\gamma_2/2}[\psi_n(x,t;\lambda)+\psi_n(x,t;\lambda')] & \text{if $\alpha=1$}
      \end{dcases} \\
      & \quad +C|x-\lambda(t+1)|^{-(\alpha_n+\beta-1)/2}|x-\lambda'(t+1)|^{-\min(\alpha,1)/2-1/2}\chi_K(x,t;\lambda,\lambda'),
    \end{aligned}
  \end{align}
  where $\gamma_2=\min(\alpha,1)+\beta-1$.
\end{lem}

\begin{proof}
  This is also proved by an adaptation of the proof of~\cite[Lemma~A.8]{Koike22b}.
\end{proof}

\begin{lem}
  \label{lem:conv_psi_n_different_family_0_sqrt}
  Let $\lambda \neq \lambda'$, $\mu>0$, $\alpha \geq 0$, and $0\leq \beta<3-\alpha_n$. Then we have
  \begin{align}
    & \int_{0}^{t^{1/2}}\int_{-\infty}^{\infty}(t-s)^{-1-\alpha}e^{-\frac{(x-y-\lambda(t-s))^2}{\mu(t-s)}}(s+1)^{-\beta/2}\psi_n(y,s;\lambda')\, dyds \\
    & \quad \leq C(t+1)^{-\alpha-(\beta-\alpha_n+1)/4}\psi_n(x,t;\lambda).
  \end{align}
\end{lem}

\begin{proof}
  This is a simple generalization of~\cite[Lemma~A.2]{Koike22b}.
\end{proof}

\begin{lem}
  \label{lem:conv_psi_n_different_family_sqrt_t_half}
  Let $\lambda \neq \lambda'$, $\mu>0$, and $\alpha \geq 0$. Then for $t\geq $, we have
  \begin{align}
    & \int_{t^{1/2}}^{t/2}\int_{-\infty}^{\infty}(t-s)^{-1/2-\alpha}e^{-\frac{(x-y-\lambda(t-s))^2}{\mu(t-s)}}(s+1)^{-\beta_{n+1}}\psi_n(y,s;\lambda')\, dyds \\
    & \quad \leq C(t+1)^{-\alpha}\psi_n(x,t;\lambda).
  \end{align}
\end{lem}

\begin{proof}
  A simple adaptation of the proof of~\cite[Lemma~A.6]{Koike22b} proves the lemma.
\end{proof}

\begin{lem}
  \label{lem:conv_exp}
  Let $\lambda \in \mathbb{R}$, $\mu>0$, and $\alpha \geq 0$. Then
  \begin{equation}
    \int_{0}^{t}e^{-\frac{t-s}{\mu}}(s+1)^{-\alpha}\psi_n(x,s;\lambda)\, ds\leq C(t+1)^{-\alpha}\psi_n(x,t;\lambda).
  \end{equation}
\end{lem}

\begin{proof}
  The lemma can be proved by slightly modifying the proof of~\cite[Lemma~3.9]{LZ97}.
\end{proof}

\subsubsection*{Acknowledgements}
This work was supported by Grant-in-Aid for JSPS Research Fellow (Grant Number 20J00882) and JSPS Grant-in-Aid for Early-Career Scientists (Grant Number 22K13938).

\bibliographystyle{amsplain}
\bibliography{kai-2021-2}

\end{document}